\documentclass[aop,MSNbibl,citesort,seceqn,dvips]{arximspdf}
\usepackage{graphicx}

\doi{10.1214/10-AOP537}
\volume{38}
\issue{6}
\pubyear{2010}
\firstpage{2258}
\lastpage{2294}

\makeatletter

\newcommand{\abs}[1]{\vert#1\vert}
\newcommand{\norm}[1]{\Vert#1\Vert}
\newcommand{\eqref}[1]{(\ref{#1})}
\newcommand{\Var}{\operatorname{Var}}
\newcommand{\Cov}{\operatorname{Cov}}
\newcommand{\iint}{\int\!\!\!\int}

\newproclaim{asm}{Assumption}
\newtheorem{theorem}{Theorem}[section]
\newtheorem{cor}[theorem]{Corollary}
\newproclaim{rem}[theorem]{Remark}
\newtheorem{prop}[theorem]{Proposition}
\newtheorem{lem}[theorem]{Lemma}
\makeatother

\begin{document}
\begin{frontmatter}

\title{Current fluctuations of a system of one-dimensional random walks
in random~environment\thanksref{T1}}
\thankstext{T1}{This work was done while both authors were visiting
Institut Mittag-Leffler (Djursholm, Sweden) for the program ``Discrete
Probability,'' and while the first author was visiting the University
of Wisconsin--Madison as a Van Vleck Visiting Assistant Professor.}
\runtitle{Current fluctuations of RWRE}

\begin{aug}
\author[A]{\fnms{Jonathon} \snm{Peterson}\ead[label=e1]{peterson@math.cornell.edu}\thanksref{t2}}
\and
\author[B]{\fnms{Timo} \snm{Sepp\"al\"ainen}\corref{}\ead[label=e2]{seppalai@math.wisc.edu}\thanksref{t3}}
\thankstext{t2}{Supported in part by NSF Grant
DMS-08-02942.}
\thankstext{t3}{Supported in part by
NSF Grant DMS-07-01091 and by the Wisconsin Alumni
Research Foundation.}
\runauthor{J. Peterson and T. Sepp\"al\"ainen}
\affiliation{Cornell University and University of Wisconsin--Madison}
\address[A]{Department of Mathematics\\
Cornell University\\
Malott Hall\\
Ithaca, New York 14850\\
USA\\
\printead{e1}} %adresu isvedimo komanda gale!
\address[B]{Department of Mathematics\\
University of Wisconsin--Madison\\
Van Vleck Hall, 480 Lincoln Dr.\\
Madison, Wisconsin 53706\\
USA\\
\printead{e2}}
\end{aug}

% HISTORY:
\received{\smonth{4} \syear{2009}}
\revised{\smonth{12} \syear{2009}}

% ABSTRACT
%
\begin{abstract}
We study the current of particles that move independently
in a common static random environment on the one-dimensional integer lattice.
A two-level fluctuation picture appears. On the central limit scale the quenched
mean of the current process converges to a Brownian motion. On a smaller
scale the current process centered at its quenched mean converges to
a mixture of Gaussian processes.
These Gaussian processes are similar to those arising from classical
random walks, but the environment makes itself felt
through an additional Brownian random shift in the spatial argument of the
limiting current process.
\end{abstract}
%
% KEYWORDS
%

\end{frontmatter}
%

%s1 ###
\section{Introduction}

We investigate the effect of a
random environment on the fluctuations of particle
current in a system of many particles.
We take the standard model of random walk in random environment
(RWRE) on the one-dimensional integer lattice, and let a large number of
particles evolve independently of each other but in a common, fixed
environment $\omega$. On the
level of the averaged (annealed) distribution
particles interact with each other through the environment.

We set the parameters of the model so that an individual particle has a positive
asymptotic speed $\mathrm{v}_P$ and satisfies a central limit theorem around
this limiting velocity
under the averaged distribution. There is also a quenched central limit theorem
that requires an environment-dependent correction $Z_n(\omega)$ to the
asymptotic value $n\mathrm{v}_P$. We scale space and time by the same
factor $n$.
We consider initial particle configurations whose distribution may
depend on
the environment, but in a manner that respects spatial shifts. Under a fixed
environment the initial occupation variables are required to be independent.

We find a two-tier fluctuation picture. On the scale ${n}^{1/2}$ the
quenched mean
of the current process behaves like a Brownian motion. In fact, up to
$o(n^{1/2})$ deviations, this quenched mean coincides with the quenched
CLT correction
$Z_n(\omega)$ multiplied by the mean density of particles.
Around its quenched mean, the current process fluctuates on the
scale $n^{1/4}$. These fluctuations are described by the same
self-similar Gaussian
processes that arise for independent particles performing classical
random walks. But the environment-determined correction $Z_n(\omega)$
appears again,
this time as an extra shift in the spatial argument of the limit
process of the current.

The broader context for this paper is the ongoing work to elucidate the patterns
of universal current fluctuations in one-dimensional driven particle systems.
A key object is the flux function $H(\mu)$ that gives the average rate
of mass
flow past a fixed point in space when the system is in a stationary
state with
mean density $\mu$.
Known rigorous results have confirmed the following delineation.
If $H $ is strictly convex or concave,
then current fluctuations have magnitude $n^{1/3}$ and limit
distributions are related to Tracy--Widom distributions from random
matrix theory.
If $H$ is linear, then the magnitude of current fluctuations
is $n^{1/4}$ and limit distributions are Gaussian.

%Thus the RWRE model confirms these predictions, for the flux
%is linear but the random environment brings some
%additional features into the picture.

The RWRE model has a linear flux. Our results show that in a sense it
confirms the prediction
stated above, but with additional features coming from the random environment.
Limit processes possess covariances that are similar to those that
arise for independent
classical random walks. However, when the environment is averaged out,
limit distributions can fail to be Gaussian.

\subsection*{Literature} A standard reference on the basic RWRE model is \cite{zRWRE}.
Further references to RWRE work follow below when we review basic results.
Earlier related results for current fluctuations of independent
particles appeared in
papers \cite{durr-gold-lebo,kSTCP} and \cite{sepp-rw}.
A central model for the study of fluctuations in the case of a concave
flux is the
asymmetric exclusion process. Key papers include \cite{bala-sepp-aom,ferr-spoh-06,joha} and \cite{quas-valko}.

Though not a system with drift, the symmetric simple exclusion process
shares some features with this class of systems with linear flux.
Namely, in the stationary process current fluctuations have magnitude
$t^{1/4}$ and fractional Brownian motion limits. This line of work
began with \cite{arratia},
with most recent contributions that give process level limits in
\cite{jara-landim-06} and \cite{peli-seth}. Fluctuations of symmetric
systems have also been
studied with disorder on the bonds \cite{jara-09,jara-landim-08}.

\subsection*{Organization of the paper} We define the model and state the
results for
the current process and its quenched mean
in Section \ref{resultsec}. Section \ref{CLTreview} reviews known
central limit results for the walk itself that we need for the proof.
Sections \ref{qmeansec}
and~\ref{currsec} prove the fluctuation theorems for the current.
An \hyperref[UIapp]{Appendix} proves a uniform integrability result for the walk that is
used in the proofs.

%s2 ###
\section{Description of the model and main results}\label{resultsec}
%In this paper we study current fluctuations in a system of independent
%one-dimensional random walks in a common random environment.
We begin with the standard RWRE model on $\mathbb Z$ with the extra feature
that we
admit infinitely many particles.
Let $\Omega:= [0,1]^{\mathbb Z}$ be the space of environments. For any
environment $\omega= \{ \omega_x \}_{x\in\mathbb Z} \in\Omega$ and
any $x\in\mathbb Z$, let
$\{X^{m,i}_\centerdot\}_{m,i}$ be a family of Markov chains with
distribution $P_\omega$ given by the following properties:
\begin{enumerate}[(1)]
\item[(1)]$\{X^{m,i}_\centerdot\}_{m \in\mathbb Z, i\in\mathbb N}$ are
independent under
the measure $P_\omega$.
\item[(2)]$P_\omega(X^{m,i}_0 = m ) = 1$, for all $m\in\mathbb Z$ and
$i\in\mathbb N$.
\item[(3)] The transition probabilities are given by
\[
P_\omega(X^{m,i}_{n+1} = x+1 | X^{m,i}_n = x ) = 1- P_\omega
(X^{m,i}_{n+1} =
x-1 | X^{m,i}_n = x ) = \omega_x.
\]
\end{enumerate}
A system of random walks in a random environment may then be
constructed by first choosing an environment $\omega$ according to a
probability distribution $P$ on $\Omega$ and then constructing the
system of random walks $\{X^{m,i}_\centerdot\}$ as described above.
The distribution $P_\omega$ of the random walks given the environment
$\omega$
is called the \textit{quenched law}.
The \textit{averaged law} $\mathbb{P}$ (also called the annealed law)
is obtained
by averaging the quenched law over all environments. That is, $\mathbb
{P}(\cdot
) := \int_{\Omega} P_\omega(\cdot) P(d\omega)$.

Often we will be considering events that only concern the behavior of a
single random walk started at location $m$, and so we will use the
notation $X^{m}_n$ in place of $X^{m,1}_n$. Moreover, if the random
walk starts at the origin, we will further abbreviate the notation by
$X_n$ in place of $X_n^0$. Expectations with respect to the measures
$P$, $P_\omega$ and $\mathbb{P}$ will be denoted by $E_P$, $E_\omega
$ and $\mathbb{E}$,
respectively, and variances with respect to the measure $P_\omega$
will be
denoted by $\Var_\omega$. Generic probabilities and expectations not
defined in the RWRE model are denoted
by $\mathbf{P}$ and $\mathbf{E}$.

For the remainder of the paper we will make the following assumptions
on the distribution $P$ of the environments.
\begin{asm}\label{UEIIDasm}
The distribution on environments is i.i.d. and uniformly elliptic.
That is,
the variables $\{\omega_x\}_{x\in\mathbb Z}$ are independent and identically
distributed under the measure $P$, and there exists a $\kappa> 0$ such
that $P(\omega_x \in[\kappa, 1-\kappa] ) = 1$.
\end{asm}
\begin{asm}\label{CLTasm}
$E_P (\rho_0^2) < 1$, where $\rho_x := \frac{1-\omega_x}{\omega_x}$.
\end{asm}

%
%The above assumptions are enough to imply a quenched central limit
%theorem with random centering. That is,
The above assumptions on the distribution $P$ on environments imply
that the RWRE are transient to $+\infty$ with strictly positive speed
$\mathrm{v}_P$ \cite{sRWRE}. That is,
%
%
%e1 ###
\begin{equation}\label{LLN}
\lim_{n\rightarrow\infty} \frac{X_n}{n} = \frac{1-E_P \rho
_0}{1+E_P \rho_0} =:
\mathrm{v}_P> 0,\qquad  \mathrm{\mathbb{P}\mbox{-}a.s.}
\end{equation}
Moreover, Assumptions \ref{UEIIDasm} and \ref{CLTasm} imply that a
quenched central limit theorem holds with a random (depending on the
environment) centering. That is, there exists an explicit function of
the environment $Z_n(\omega)$ and a constant $\sigma_1 > 0$ such that for
$P\mathrm{\mbox{-}a.e.}$ environment $\omega$,
\[
\lim_{n\rightarrow\infty} P_\omega\biggl( \frac{X_n - n\mathrm
{v}_P+Z_n(\omega)}{\sigma_1 \sqrt{n}}
\leq x \biggr) = \Phi(x)\qquad  \forall x\in\mathbb R,
\]
where $\Phi$ is the standard normal distribution function.
The environment-dependent centering in the above quenched central limit
theorem cannot be replaced by a deterministic centering since it is
known that there exists a constant $\sigma_2 > 0$\vspace*{1pt} such that the process
$t\mapsto\frac{Z_{nt}(\omega)}{\sigma_2 \sqrt{n}}$ converges
weakly to a
standard Brownian motion.\vspace*{1pt} Definitions of $\sigma_1,\sigma_2$ and
$Z_n(\omega)$ are
provided in Section \ref{CLTreview} where we give a more detailed
review of the known limit distribution results for RWRE under
Assumptions \ref{UEIIDasm} and~\ref{CLTasm}.

In this paper we will be concerned with a system of RWRE in a common
environment with a finite (random) number of walks started at each site
$x\in\mathbb Z$. Let $\eta_0(x)$ be the number of walks started from
$x \in\mathbb Z$.
%, and
%be the number of the random walks at location $x$ at time $t$.
We will allow the law of the initial configurations to depend on the
environment (in a measurable way). Let $\theta$ be the shift operator
on environments defined by $(\theta^x\omega)_y = \omega_{x+y}$. We
will assume
that our initial configurations are stationary in the following sense.
\begin{asm}\label{ICasm}
The distribution of $\eta_0$ is such that $\omega\mapsto P_\omega
(\eta_0(0) =
k)$ is a measurable function of $\omega$ for any $k\in\mathbb N$, and
the law of
$\eta_0$ respects the shifts of the environment: $P_\omega( \eta
_0(x) = k )
= P_{\theta^x\omega}( \eta_0(0) = k )$. Also, given the environment~$\omega$,
the $\{ \eta_0(x) \}$ are independent and independent of the paths of
the random walks.
\end{asm}

We will also need the following moment assumptions.
\begin{asm}\label{ICMasm}
For some $\varepsilon>0$,
%
%
%e2 ###
\begin{equation}
E_P [ E_\omega(\eta_0(x))^{2+\varepsilon} + \Var_\omega(\eta
_0(x))^{2+\varepsilon} ]<\infty.
\label{vpmomass1}
\end{equation}
\end{asm}

To simplify notation, we will let $\bar\mu(\omega) := E_\omega[
\eta_0(0)]$. Note
that Assumption \ref{ICasm} implies that $E_\omega[ \eta_0(m) ] =
\bar\mu
(\theta^m\omega)$.
Let $\mu:= E_P[ \bar\mu(\omega)] = \mathbb{E}\eta_0(0)$ be the
average density of
the initial configuration of particles, and let $\sigma_0^2=E_P[
\Var_\omega(\eta_0(x)) ]$.

The law of large numbers \eqref{LLN} implies that each random walk
moves with asymptotic speed $\mathrm{v}_P$.
The main object of study in this paper is the following two-parameter
process. For $t\geq0$ and $r\in\mathbb R$, let
%
%
%e3 ###
\begin{eqnarray}\label{defYn}
Y_n(t,r) &=& \sum_{m> 0} \sum_{k=1}^{\eta_0(m)} \mathbf{1}{ \bigl\{
X^{m,k}_{nt} \leq nt \mathrm{v}_P+r\sqrt n  \bigr\} }
\nonumber
\\[-8pt]
\\[-8pt]
\nonumber
&&{}- \sum_{m\leq0} \sum_{k=1}^{\eta_0(m)} \mathbf{1}{ \bigl\{ X^{m,k}_{nt}
> nt\mathrm{v}_P+ r\sqrt n   \bigr\} } .
\end{eqnarray}
\begin{figure}

\includegraphics{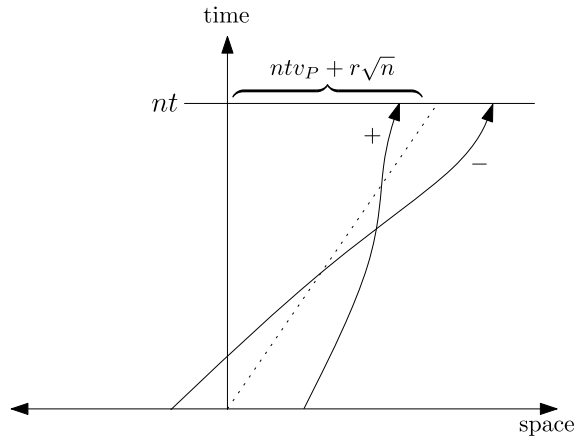}

\caption{A visual representation of the process $Y_n(t,r)$ which is the net
(negative) current seen by an observer starting at the origin at time 0
and ending at $nt\mathrm{v}_P+ r\sqrt{n}$ at time $nt$. Particles
crossing from
right to left contribute positively and those crossing from left to
right contribute negatively.}\label{fig1}\vspace*{-3pt}
\end{figure}

A visual description of the process $Y_n(t,r)$ is given in Figure 1.
$Y_n(t,r)$ is similar to what was called the space--time current process
in \cite{kSTCP}
and studied in a constant environment (i.e., particles performing
independent classical random walks).
We altered the definition because the limit process of this version has
a more
natural description. The process studied earlier in \cite{kSTCP} equals\vspace*{-2pt}
\begin{eqnarray}\label{oldcurrent}
Y_n(t,r) - Y_n(0,r) &=& \sum_{m> r\sqrt n} \sum_{k=1}^{\eta_0(m)} \mathbf{1}{\bigl \{
X^{m,k}_{nt} \leq nt\mathrm{v}_P+r\sqrt n  \bigr\} }
\nonumber
\\[-10pt]
\\[-10pt]
\nonumber
&&{}- \sum_{m\leq
r\sqrt n} \sum_{k=1}^{\eta_0(m)} \mathbf{1}{ \bigl\{ X^{m,k}_{nt} >
nt\mathrm{v}_P+ r\sqrt n   \bigr\} } .
\end{eqnarray}
This process $Y_n(\cdot,r) - Y_n(0,r)$ is the net right-to-left
particle current seen by an observer who starts at $r\sqrt{n}$ and
moves with deterministic speed $\mathrm{v}_P$.
Adapting the proof of \cite{kSTCP} to our definition of $Y_n(t,r)$
gives this theorem:
\begin{theorem}[(Kumar \cite{kSTCP})]
Assume that the environment is nonrandom.\vspace*{1pt} That is, there exists a $p\in
(0,1)$ such that $P(\omega_x = p, \forall x\in\mathbb Z) = 1$. Let
$\mathbb{E}(\eta_0) =
\mu$ and ${\rm\mathbb{V}ar}(\eta_0) = \sigma_0^2$, and assume that
$\mathbb{E}(\eta_0^{12}) < \infty$. Then, the process $n^{-1/4}(
Y_n(\cdot,
\cdot) - \mathbb{E}Y_n(\cdot,\cdot) )$ converges in distribution on the
$D$-space of two-parameter cadlag processes.
The limit is the mean zero Gaussian process $V^0(\cdot,\cdot)$ with covariance
%
%
%e4 ###
\begin{equation}\label{CRWcov}
\mathbf{E}[ V^0(s,q)V^0(t,r) ] = \Gamma((s,q),(t,r)),
\end{equation}
where the covariance function $\Gamma$ is defined below in \eqref{Gammadef}.
%Moreover, if the initial occupation variables $\{\eta_0(x)\}_{x\in\Z}$
%are \iid Poisson$(\rho)$ random variables, then the process $n^{-1/4}
%Y_n(\cdot,0)$ converges to a fractional Brownian motion with Hurst
%parameter $1/4$.
\end{theorem}

The theorem above uses the higher moment assumption $\mathbb{E}(\eta
_0^{12}) <
\infty$
for process-level
tightness. We have not proved such tightness, hence, we get by with the
moments assumed
in \eqref{vpmomass1}. We turn to discuss the results in the random environment.\vadjust{\goodbreak}
%In \cite{kSTCP}, the space-time current process $Y_n(t,r)$ was studied
%in the case of independent random walks in a deterministic homogeneous
%environment (not in a random environment). In that case, it was shown
%that the fluctuations of $Y_n(t,r)$ were of the order $n^{1/4}$.
%Moreover, $n^{-1/4}Y_n(\cdot, \cdot)$ converges in distribution to a
%mean zero Gaussian process with an explicit (although somewhat
%complicated) covariance function. If the initial configuration $
%fixed $r$ the process $n^{-1/4} Y_n(\cdot,r)$ converges to fractional
%Brownian motion with Hurst parameter $1/4$.

The random environment adds a new layer
of fluctuations to the current.
%substantially changes the fluctuations of the current.
%In fact, $Y_n(t,r)$ is approximately equal to the random centering
%term $Z_{nt}(\w)$ from the quenched CLT (which is of order $\sqrt{n}$).
These larger fluctuations are of order $\sqrt{n}$ and depend only on
the environment.
This is summarized by our first main result. The process $Z_{nt}(\omega)$
in the statement below
is the correction required in the quenched central limit theorem of the
walk, defined in \eqref{Zdef}
in Section \ref{CLTreview}.
\begin{theorem} \label{QMCurrent}
For any $\varepsilon>0$, $0<R, T<\infty$,
%
%
%e5 ###
\begin{equation}\qquad\lim_{n\to\infty} P\Bigl(   \sup_{t\in[0,T],
r\in[-R,R]}
\bigl\vert E_\omega Y_n(t,r) - \mu r \sqrt{n} - \mu Z_{nt}(\omega)
\bigr\vert
\ge\varepsilon\sqrt n  \Bigr) =0.
\label{limYZ}
\end{equation}
Moreover, since $\{ n^{-1/2} Z_{nt}(\omega)\dvtx  t\in\mathbb R_+\}$
converges weakly
to $\{ \sigma_2 W(t)\dvtx  t\in\mathbb R_+ \}$, where $W(\cdot)$ is a
standard Brownian
motion, then the two-parameter process $\{ n^{-1/2} E_\omega Y_n(t,r)\dvtx
t\in
\mathbb R_+,   r\in\mathbb R\}$ converges weakly to $\{ \mu\sigma_2
W(t) + \mu r\dvtx
t\in\mathbb R_+,   r\in\mathbb R\}$.
%Also, for any fixed $r\in\R$, $n^{-1/2} (Y_n(\cdot, r) - E_\w Y_n(
%zero.
%In particular, this implies that for any fixed $r$, $\frac{1}{
%where $W(\cdot)$ is a standard Brownian motion.
\end{theorem}

To see the next order of fluctuations, we center the current at its
quenched mean. Define
%
%
%e6 ###
\begin{eqnarray}\label{vpdefVn}
V_n(t,r) &=&Y_n(t,r) -E_\omega Y_n(t,r) \nonumber\\
&=& \sum_{m> 0}
\Biggl(  \sum_{k=1}^{\eta_0(m)} \mathbf{1}{ \bigl\{ X^{m,k}_{nt} \leq
nt\mathrm{v}_P+r\sqrt{n}  \bigr\} }\nonumber\\
&&\qquad {} - E_\omega(\eta_0(m))P_\omega\bigl\{
X^m_{nt} \leq nt\mathrm{v}_P+r\sqrt{n}  \bigr\}
\Biggr)\\
&&{}  - \sum_{m\leq0} \Biggl(  \sum_{k=1}^{\eta
_0(m)}\mathbf{1}{ \bigl\{ X^{m,k}_{nt} > nt\mathrm{v}_P+r\sqrt{n}  \bigr\}
}\nonumber\\
&&\qquad\hspace*{14pt} {}- E_\omega(\eta_0(m))P_\omega\bigl\{ X^{m}_{nt} >
nt\mathrm{v}_P+r\sqrt{n}  \bigr\} \Biggr).\nonumber
\end{eqnarray}
The fluctuations of $V_n(t,r)$ are of order $n^{1/4}$ and the same as
the current fluctuations in a deterministic environment, up to a random
shift coming from the environment.
We need to introduce some notation. For any $\alpha>0$, let $\phi
_{\alpha
^2}(\cdot)$ and $\Phi_{\alpha^2}(\cdot)$ be the density and distribution
function, respectively, for a Gaussian distribution with mean zero and
variance $\alpha^2$.
Also, let
%
%
%e7 ###
\begin{eqnarray}\label{Psidef}
\Psi_{\alpha^2}(x) &:=& \alpha^2 \phi_{\alpha^2}(x) - x \Phi
_{\alpha^2}(-x)\quad   \mbox{and}
\nonumber
\\[-8pt]
\\[-8pt]
\nonumber
 \Psi_0(x) &:=& \lim_{\alpha\rightarrow0} \Psi_{\alpha
^2}(x) = x^-.
\end{eqnarray}
Then, for any $(s,q),(t,r) \in\mathbb R_+ \times\mathbb R$ define the
covariance function
%
%
%e8 ###
\begin{eqnarray}\label{Gammadef}
\Gamma((s,q),(t,r)) &:=& \mu\bigl( \Psi_{\sigma_1^2(s+t)}(q-r) -
\Psi_{\sigma
_1^2|s-t|}(q-r) \bigr)
\nonumber
\\[-8pt]
\\[-8pt]
\nonumber
&&{} + \sigma_0^2 \bigl( \Psi_{\sigma_1^2 s}(-q) + \Psi_{\sigma
_1^2 t} (r) - \Psi
_{\sigma_1^2(s+t)}(r-q) \bigr),
\end{eqnarray}
where $\sigma_1$ is the scaling factor in the quenched central limit theorem
[see \eqref{clt1} in Section \ref{CLTreview} for a formula].
Given the above definitions, let $(V,Z)=( V(t,r), Z(t)\dvtx  t\in\mathbb
R_+, r\in\mathbb R
)$ be the process whose joint
distribution is defined as follows:
\begin{longlist}[(ii)]
\item[(i)] Marginally, $Z(\cdot)=\sigma_2 W(\cdot)$ for a standard
Brownian motion $W(\cdot)$,
and $\sigma_2$ is the scaling factor in the central limit theorem of
the correction
$Z_{nt}(\omega)$
[see \eqref{clt2} in Section \ref{CLTreview} for a formula].
\item[(ii)] Conditionally on the path $Z(\cdot)\in C(\mathbb
R_+,\mathbb R)$, $V$ is
the mean zero Gaussian process indexed
by $\mathbb R_+\times\mathbb R$ with covariance
%
%
%e9 ###
\begin{eqnarray}\label{vpcov}
&&\mathbf{E}[V(s,q)V(t,r) \vert  Z(\cdot)]
\nonumber
\\[-8pt]
\\[-8pt]
\nonumber
&&\qquad=\Gamma
\bigl(\bigl(s,q+Z(s)\bigr),\bigl(t,r+Z(t)\bigr)\bigr)\qquad
 \mbox{for $(s,q), (t,r)\in\mathbb R_+\times\mathbb R$.}
\end{eqnarray}
\end{longlist}
An equivalent way to say this is to first take independent $(V^0,Z)$
with $Z$ as above
and
$V^0=\{V^0(t,r)\dvtx  (t,r)\in\mathbb R_+\times\mathbb R\}$ the mean zero
Gaussian process
with
covariance $ \Gamma((s,q),(t,r))$ from \eqref{Gammadef},
and then
define $V(t,r)=V^0(t,r+Z(t))$.

The next theorem gives joint convergence of the centered current
process and the
environment-dependent shift.

\begin{theorem}
Under the averaged probability\vspace*{1.5pt} $\mathbb{P}$, as $n\to\infty$,
the finite-dimensional distributions of the joint process
$\{( n^{-1/4}V_n(t,r), n^{-1/2} Z_{nt}(\omega) )\dvtx  t\in\mathbb
R_+, r\in\mathbb R
\}$
converge to those of the process
$(V,Z)$.
\label{findimthm}
\end{theorem}

Our proof shows additionally that
\begin{eqnarray*}
&&\lim_{n\to\infty}
E_P\Biggl\vert  E_\omega\exp\Biggl\{in^{-1/4} \sum_{k=1}^N \alpha_k
V_n(t_k,r_k) \Biggr\}\\
 &&\hspace*{41pt}\qquad{} -   \mathbf{E}\exp\Biggl\{i \sum_{i=1}^N \alpha_k V(t_k,r_k)
\Biggr\}
\Biggr\vert=0
\end{eqnarray*}
for any choice of time--space points $(t_1,r_1),\ldots,(t_N,r_N)\in
\mathbb R
_+\times\mathbb R$\vspace*{1pt} and
$\alpha_1,\ldots,\break \alpha_N \in\mathbb R$. [See \eqref{vplim7} below.]
This falls short of a quenched limit for $n^{-1/4}V_n$ (a~limit
for a fixed $\omega$), but it does imply that if a quenched limit exists,
the limit
process is the one that
we describe. We suspect, however, that no quenched limit exists since
the techniques of this paper can be used to show that the quenched
covariances of the process $n^{-1/4}V_n(\cdot,\cdot)$ do not converge $P\mathrm{\mbox{-}a.s.}$

The mean zero Gaussian process $\{u(t,r)\dvtx  t\in\mathbb R_{+},r\in
\mathbb R\}$ with
covariance
$ \mathbf{E}[u(s,q)u(t,r)]= \Gamma((s,q),(t,r))$ from
\eqref{Gammadef}
can be represented as the sum of two integrals:
%
%
%e10 ###
\begin{eqnarray}\label{udef}
u(t,r)&=& \sqrt{\mu} \iint_{[0,t]\times\mathbb R}
\phi_{\sigma_1^2(t-s)}(r-x) \,dW(s,x)
\nonumber
\\[-8pt]
\\[-8pt]
\nonumber
&&{}+\sigma_0\int_{\mathbb R}
\phi_{\sigma_1^2t}(r-x)B(x) \,dx,
\end{eqnarray}
where
$W$ is a two-parameter
Brownian motion on $\mathbb R_+ \times\mathbb R$ (Brownian sheet) and~$B$ an independent
two-sided one-parameter Browian motion on $\mathbb R$.
The process $u(t,r)$ is also a weak solution of the
stochastic heat equation with initial data given by Brownian motion
\cite{walsh}:
%
%
%e11 ###
\begin{equation}
\quad u_t = \frac{\sigma_1^2}2 u_{rr} + \sqrt{\mu}   \dot W ,\qquad
u(0,r)=\sigma_0 B(r),\qquad  (t,r)\in\mathbb R_{+}\times\mathbb R.
\label{stheateq2}
\end{equation}
This type of process we obtain if we define
$u(t,r)=V(t,r-Z(t))$ by regarding
the random path $-Z(\cdot)$ as the new spatial origin.

We next remark on the distribution of the limiting process $V(t,r)$ in
a couple of special cases. First we consider the case when $\sigma_0=0$
(this includes the case of deterministic initial configurations). If
$\sigma
_0=0$, then \eqref{Gammadef} and \eqref{vpcov} imply that, for any
fixed $t\geq0$, the one-parameter process $V(t,\cdot)$ has conditional
covariance
\begin{eqnarray*}
E[V(t,q)V(t,r) \vert  Z(\cdot)]&=&\Gamma
\bigl(\bigl(t,q+Z(t)\bigr),\bigl(t,r+Z(t)\bigr)
\bigr) \\
&=& \mu\bigl(\Psi_{2\sigma_1^2 t}(q-r) - \Psi_{0}(q-r)\bigr).% =
\end{eqnarray*}
%
%$E[V(t,q)V(t,r) \vert  Z(\cdot)]=\Gamma((t,q+Z(t)),(t,r+Z(t))
%) = \dc\Psi_{2\s_1^2 t}(|q-r|)$.
In particular, the covariances of $V(t,\cdot)$ do not depend on the
process $Z(\cdot)$ and are the same as in the classical random walk case.
\begin{cor}
If $\sigma_0 = 0$, then for any fixed $t\geq0$ the (averaged)
finite-dimensional distributions of the one parameter process $\{ n^{-1/4}
V_n(t,r) \dvtx   r\in\mathbb R\}$ converge to those of the one parameter mean
zero Gaussian process $V^0(t,\cdot)$ with covariances given by \eqref
{CRWcov} with $s=t$.
\end{cor}

A second special case worth considering is when $\mu= \sigma_0^2$.
%In this case \eqref{Gammadef} simplifies to
%$\Gamma((s,q),(t,r)) = \dc( \Psi_{\s_1^2 s}(-q) +
In the case of classical random walks, $\mu= \sigma_0^2$ implies that
\[
\mathbf{E}[ V^0(s,0)V^0(t,0) ] = \frac{\mu\sigma_1}{\sqrt{2\pi}}
\bigl( \sqrt{s} + \sqrt
{t} - \sqrt{|s-t|}\bigr),
\]
so that $V^0(\cdot,0)$ is a fractional Brownian motion with Hurst
parameter $1/4$.
For RWRE, $\mu= \sigma_0^2$ implies that
%
%
%e12 ###
\begin{eqnarray}\label{condcov}
&&\mathbf{E}[V(s,0)V(t,0) \vert  Z(\cdot)]
\nonumber
\\[-8pt]
\\[-8pt]
\nonumber
%=\Gamma((s,Z(s)),(t,Z(t)))
&&\qquad= \mu\bigl( \Psi_{\sigma_1^2 s}(-Z(s)) + \Psi_{\sigma_1^2
t}(Z(t)) - \Psi_{\sigma
_1^2|s-t|}\bigl(Z(t)-Z(s)\bigr) \bigr).
\end{eqnarray}
Since the right-hand side of \eqref{condcov} is a nonconstant random
variable, the marginal distribution of $V(t,0)$ is non-Gaussian. Taking
expectations of \eqref{condcov} with respect to $Z(\cdot)$ gives that
%
%
%e13 ###
\begin{equation}\label{specialcasecov}
\mathbf{E}[V(s,0)V(t,0)] = \frac{\mu\sqrt{\sigma_1^2 + \sigma
_2^2}}{\sqrt{2\pi}} \bigl(
\sqrt{s} + \sqrt{t} - \sqrt{|s-t|}\bigr).
\end{equation}
Thus, we have the following.
\begin{cor}
If $\mu= \sigma_0^2$, then the process $V(\cdot,0)$ has covariances like
that of a fractional Brownian motion, but is not a Gaussian process.
\end{cor}
\begin{rem}
The condition that $\mu= \sigma_0^2$ is important because it includes the
case when the configuration of particles is stationary under the
dynamics of the random walks. For classical random walks, the
stationary distribution on configurations of particles is when the
$\eta
_0(x)$ are i.i.d. $\operatorname{Poisson}(\mu)$ random variables. Consider now the case
where, given $\omega$, the $\eta_0(x)$ are independent and
%
%
%e14 ###
\begin{equation}\label{fdef}
\qquad\quad \eta_0(x) \sim\operatorname{Poisson}(\mu f(\theta^x\omega))\qquad
 \mbox{where }
f(\omega) =
\frac{\mathrm{v}_P}{\omega_0}\Biggl( 1 + \sum_{i=1}^{\infty} \prod
_{j=1}^i \rho_{j}
\Biggr).
\end{equation}
It was shown in \cite{psHydro} that, given $\omega$, the above distribution
on the configuration of particles is stationary under the dynamics of
the random walks. Note that in this case, $E_\omega\eta_0(0) = \Var
_\omega\eta
_0(0) = \mu f(\omega)$. Moreover, Assumptions \ref{UEIIDasm} and \ref
{CLTasm} imply that $E_P \rho_0^{2+\varepsilon} < 1$ for some
$\varepsilon>0$, and, thus,
it can be shown that $E_P f(\omega)^{2+\varepsilon} < \infty$. Therefore,
Assumptions \ref{ICasm} and \ref{ICMasm} are fulfilled in this
special case.
\end{rem}

It is intuitively evident but not a corollary of our theorem that if
the environment-dependent
shift is introduced in the current process itself,
the random shift $Z$
disappears from the limit process $V$. For the sake of completeness, we
state this result too.
For $(t,r)\in\mathbb R_+\times\mathbb R$ define
%
%
%e15 ###
\begin{eqnarray}\label{defYnq}
Y_n^{(q)}(t,r) &=& \sum_{m> 0} \sum_{k=1}^{\eta_0(m)}
\mathbf{1}{ \bigl\{ X^{m,k}_{nt} \leq nt\mathrm{v}_P-Z_{nt}(\omega
)+r\sqrt n  \bigr\} }
\nonumber
\\[-8pt]
\\[-8pt]
\nonumber
&&{}- \sum_{m\leq0} \sum_{k=1}^{\eta_0(m)} \mathbf{1}{ \bigl\{ X^{m,k}_{nt}
> nt\mathrm{v}_P -Z_{nt}(\omega)+ r\sqrt n   \bigr\} }
\end{eqnarray}
and its centered version
\[
V_n^{(q)}(t,r) = Y_n^{(q)}(t,r) - E_\omega Y_n^{(q)}(t,r).
\]
The process $V_n^{(q)}$ has the same limit as classical
random walks.
As above, let $V^0=\{V^0(t,r)\dvtx  (t,r)\in\mathbb R_+\times\mathbb R\}$
be the mean zero
Gaussian process
with
covariance~\eqref{CRWcov}.

\begin{theorem}
Under the averaged probability $\mathbb{P}$, as $n\to\infty$,
the finite-dimensional distributions of the joint process
$\{( n^{-1/4}V_n^{(q)}(t,r), n^{-1/2} Z_{nt}(\omega) )\dvtx\break  t\in
\mathbb R_+, r\in
\mathbb R\}$
converge to those of the process
$(V^0,Z)$ where $V^0$ and $Z$ are independent.
%Moreover, for any fixed $(t,r) \in\R_+ \times\R$,
% \lim_{n\ra\infty} E_P[ ( n^{-1/4} E_\w Y_n^{(q)}(t,r)
%)^2 ] = 0.
%Thus, under the averaged probability $\P$, the finite dimensional
%distributions of $Y_n^{(q)}(\cdot,\cdot)$ converge to those of $V^0(
\label{findimthmVq}
\end{theorem}

It can be shown, using the techniques of this paper, that $n^{-1/2}
E_\omega Y_n^{(q)}(t,r)$ converges to zero in probability for any fixed $t$
and $r$. We suspect that the fluctuations of $E_\omega Y_n^{(q)}(t,r)$ are
at most of order $n^{-1/4}$, but at this point we have no result.
%are unable to determine the limiting distribution of $n^{-1/4} E_\w
%Y_n^{(q)}(t,r)$. In particular, it may be that this converges to zero
%in probability, in which case the above theorem would imply that the
%finite dimensional distributions of $Y_n^{(q)}(\cdot,\cdot)$ are the
%same as in the case of classical random walks.

%At this point we are unable to identify the limiting distribution of
%$E_\w Y_n^{(q)}(t,r)$. Using the techniques of this paper, it is
%possible to show that the fluctuations of $E_\w Y_n^{(q)}(t,r)$ are of
%the order $n^{1/4}$. In fact,
%& \lim_{n\ra\infty} E_P [ ( n^{-1/4} E_\w Y_n^{(q)}(t,r)
%)^2 ] \\
%& = \dc^2 ( 2 \int_0^{\infty} \Phi_{\s_1^2 t}(-x)^2 dx - 2
%( \int_0^\infty\Phi_{\s_1^2 t}(-x) dx )^2 + 4 \int_0^
%)\\
%& = \s_1^2 t( 1-\frac{2}{\pi} ) + \s_1 \sqrt{t} (
%Currently, we are not able to compute the limit of $E_P [
%n^{-1/4} E_\w Y_n^{(q)}(t,r) ]$, nor are we able to calculate
%the limiting distribution. However, we suspect that the limiting mean
%is zero and that the limiting distributions are Gaussian.

%s3 ###
\section{Review of CLT for RWRE}\label{CLTreview}

In this section we review some of the limiting distribution results for
one-dimensional RWRE implied by Assumptions \ref{UEIIDasm} and \ref
{CLTasm}. Before stating a theorem which summarizes what is known, we
introduce some notation.
Let $T_x := \inf\{ n\geq0\dvtx  X_n = x \}$ be the hitting time of the site
$x\in\mathbb Z$ of a~RWRE started at the origin, and for $x\in\mathbb
Z$ let
%
%
%e16 ###
\begin{equation}\label{hdef}
h(x,\omega) :=
\cases{
\mathrm{v}_P\displaystyle\sum_{i=0}^{x-1} ( E_{\theta^i \omega} T_1 - \mathbb
{E}T_1 ), & \quad $x \geq1$, \vspace*{2pt}\cr
0, & \quad  $x= 0$, \vspace*{2pt}\cr
- \mathrm{v}_P\displaystyle\sum_{i=x}^{-1} ( E_{\theta^i \omega} T_1 - \mathbb
{E}T_1 ), & \quad $x \leq-1$.}
%
%, \mbox{and}
%Z_n(\w) := h(\fl{n\vp},\w).
\end{equation}
Define also
%
%
%e17 ###
\begin{equation}Z_{nt}(\omega) := h(\lfloor nt\mathrm{v}_P \rfloor
,\omega). \label{Zdef}
\end{equation}
%
%The above assumptions imply the following theorem which describes the
%behavior of a single RWRE starting from the origin.
%
\begin{theorem}[(\cite{gQCLT,kksStable,pThesis,zRWRE})]\label{QCLTthm}
Let Assumptions \textup{\ref{UEIIDasm}} and \textup{\ref{CLTasm}} hold. Then, the
following hold:
\begin{enumerate}[(1)]
%That is, \begin{equation}\label{LLN}
% \lim_{n\ra\infty} \frac{X_n}{n} = \frac{1-E_P \rho_0}{1+E_P \rho_0}
%=: \vp> 0,  \P-a.s.
%
\item[(1)] The RWRE satisfies a quenched functional central limit theorem
with a random (depending on the environment) centering.
For $n\in\mathbb N$ and $t\geq0$, let
%
%
%e18 ###
\begin{equation}
B^n(t) := \frac{X_{nt} - nt\mathrm{v}_P+ Z_{nt}(\omega)}{\sigma_1
\sqrt{n}}\qquad
\mbox{where } \sigma_1^2 := \mathrm{v}_P^3 E_P( \Var_\omega T_1 ).
\label{clt1}
\end{equation}
%
%where $\s_1^2 := \vp^3 E_P( \Var_\w T_1 )$.
Then, for $P$-a.e. environment $\omega$, under the quenched measure
$P_\omega$,
$B^n(\cdot)$ converges weakly to standard Brownian motion as
$n\rightarrow\infty$.
\item[(2)] Let
%
%
%e19 ###
\begin{equation}
\zeta^n(t) := \frac{Z_{nt}(\omega)}{\sigma_2 \sqrt{n}}\qquad
\mbox{where } \sigma
_2^2 := \mathrm{v}_P^2 \Var( E_\omega T_1).
\label{clt2}
\end{equation}
Then, under the measure $P$ on environments, $\zeta^n(\cdot)$ converges
weakly to standard Brownian motion as $n\rightarrow\infty$.
%Brownian motion as $n\ra\infty$.
%
\item[(3)] The RWRE satisfies an averaged functional central limit theorem. Let
\[
\mathbb{B}^n(t) := \frac{X_{nt} - nt\mathrm{v}_P+ Z_{nt}(\omega
)}{\sigma\sqrt{n}}\qquad
 \mbox{where } \sigma^2 = \sigma_1^2 + \sigma_2^2.
\]
Then, under the averaged measure $\mathbb{P}$, $\mathbb{B}^n(\cdot)$
converges
weakly to standard Brownian motion.
\end{enumerate}
\end{theorem}
\begin{rem}
The conclusions of Theorem \ref{QCLTthm} still may hold if the law on
environments is not uniformly elliptic or i.i.d. but satisfies certain
mixing properties \cite{gQCLT,kksStable,mrzStable,pThesis,zRWRE}.
However, if the environment is i.i.d., the requirement that $E_P \rho
_0^2 < 1$ in Assumption \ref{CLTasm} cannot be relaxed in order for
Theorem \ref{QCLTthm} to hold \cite{kksStable,p1LSL2,pzSL1}.
%If $E_P \log\rho_0 < 0$ and $E_P \rho^s = 1$ for some $s\in(0,2)$
%(implying $E_P \rho^2 > 1$), then there does not exist a quenched
%limiting distribution \cite{pzSL1}, and an averaged limiting
%distribution holds with sub-diffusive scaling and a non-Gaussian limit
%theorem still holds but with scaling different than $\sqrt{n}$
\end{rem}

Let $B_\centerdot$ denote a standard Brownian motion with distribution
$\mathbf{P}$.
The quenched functional central limit theorem implies that, $P\mathrm{\mbox{-}a.s.}$,
for any $s,t\geq0$ and $x,y\in\mathbb R$,
%
%
%e20 ###
\begin{eqnarray}\label{QCLT}
&&\lim_{n\rightarrow\infty} P_\omega\biggl( \frac{X_{ns} - ns\mathrm
{v}_P+Z_{ns}(\omega)}{\sigma_1\sqrt
{n}} \leq x ,  \frac{X_{nt} - nt\mathrm{v}_P+Z_{nt}(\omega)}{\sigma
_1\sqrt{n}} \leq y
\biggr)
\nonumber
\\[-8pt]
\\[-8pt]
\nonumber
&&\qquad = \mathbf{P}( B_s \leq x,   B_t \leq y ),
\end{eqnarray}
where $B_\centerdot$ is a standard Brownian motion.
Moreover, for fixed $s,t > 0$, the convergence in \eqref{QCLT} is
uniform in $x$ and $y$.
In \cite{zRWRE}, only an averaged central limit theorem is proved.
%However, since $\mathbb{B}^n(t) = \frac{\s_1}{\s}B^n(t) + \frac{
%Z_{nt}(\w)}{\s\sqrt{n}}$,
However, since $\mathbb{B}^n(t) = \frac{\sigma_1}{\sigma}B^n(t) +
\frac{\sigma_2}{\sigma
}\zeta^n(t)$,
the averaged functional central limit theorem can be derived from the
previous two parts of Theorem \ref{QCLTthm}. Indeed, it follows
immediately that the finite-dimensional distributions of $\mathbb
{B}^n(t)$ converge to those of a Brownian motion (as in \cite{zRWRE},
this uses that convergences of terms like \eqref{QCLT} hold uniformly
in $x$ and $y$). Thus, it only remains to show that~$\mathbb
{B}^n(\cdot
)$ is tight, but this is not too difficult.

The random centering $nt\mathrm{v}_P- Z_{nt}(\omega)$ in the quenched
CLT is more
convenient than centering by the quenched mean $E_\omega X_{\lfloor nt
\rfloor}$.
Both centerings are essentially the same in the sense that they do not
differ on the scale of $\sqrt{n}$:
%
%
%e21 ###
\begin{equation}\label{limqmeanZ}
\lim_{n\rightarrow\infty} P\Bigl( \sup_{k\leq n} | E_\omega X_k -
k\mathrm{v}_P+ Z_k(\omega) |
\geq\varepsilon\sqrt{n} \Bigr) = 0\qquad   \forall\varepsilon>0.
\end{equation}
But $Z_{nt}(\omega) = h( \lfloor nt\mathrm{v}_P \rfloor,\omega)$ is
convenient because it is
defined in terms of partial sums of the
%sum of the stationary, ergodic sequence of
random variables $E_{\theta^i\omega} T_1$ for which there is an explicit
formula in terms of the environment $\omega$ (see \cite{pThesis} or
\cite{zRWRE}).
We note the following lemma due to Goldsheid \cite{gQCLT} which we will
use in several places in
the remainder of the paper.
\begin{lem}\label{Goldsheid}
Let Assumptions \textup{\ref{UEIIDasm}} and \textup{\ref{CLTasm}} hold. Then there exists
an $\eta>0$ and a~constant $C<\infty$ such that
%
%
%e22 ###
\begin{equation}\label{hLpbound}
E_P \Bigl[ \sup_{1\leq k \leq n} |h(k,\omega)|^{2+2\eta} \Bigr]
\leq C
n^{1+\eta}\qquad   \forall n\in\mathbb N.
\end{equation}
\end{lem}

We conclude\vspace*{1pt} this section by stating a new result on the uniform
integrability (under the averaged measure) of $n^{-1/2} (X_n - n\mathrm{v}_P)$.
\begin{prop}\label{UIprop}
Let $\sigma_1^2$ and $\sigma_2^2$ be defined as in Theorem \textup{\ref
{QCLTthm}}. Then,
%
%
%e23 ###
\begin{equation}\label{UIclaim}
\lim_{n\rightarrow\infty} \frac{1}{n} \mathbb{E}( X_n -
n\mathrm{v}_P)^2 = \sigma_1^2 +
\sigma_2^2.
\end{equation}
Moreover, there exists a constant $C<\infty$ such that
%
%
%e24 ###
\begin{equation}\label{supUIclaim}
\mathbb{E}\Bigl[\sup_{k\leq n} (X_k - k \mathrm{v}_P)^2 \Bigr]
\leq C n.
\end{equation}
\end{prop}

The proof of Proposition \ref{UIprop} is given in \hyperref[UIapp]{Appendix}.
It should be noted that while the statement \eqref{limqmeanZ} does not
appear anywhere in the literature (at least that we know of),
it is included in the proof of Proposition \ref{UIprop}.
%it follows easily from the same methods that are used in the proof of
%Proposition \ref{UIprop}.

%s4 ###
\section{Fluctuations of the quenched mean of the current}\label{qmeansec}
In this section we prove Theorem \ref{QMCurrent} for the quenched mean
of $Y_n(t,r)$.
%For ease of notation we will let $W_n(t,r) := E_\w Y_n(t,r)$.
Introduce the notation
\begin{eqnarray}\label{Wndef}
W_n(t,r) &:=& E_\omega Y_n(t,r) - \mu r\sqrt{n} \nonumber\\
&=& \sum_{m> 0} E_\omega[ \eta_0(m) ] P_\omega\bigl( X^m_{nt} \leq
nt\mathrm{v}_P+ r\sqrt{n}
\bigr)\\
&&{} - \sum_{m\leq0}E_\omega[ \eta_0(m) ] P_\omega\bigl( X^m_{nt} >
nt\mathrm{v}_P+ r\sqrt{n}
\bigr) - \mu r\sqrt{n}.\nonumber
\end{eqnarray}
%
% W_n(t,r) &:= E_\w Y_n(t,r) - \dc r\sqrt{n} \\
%&= \sum_{m> 0} \edc(\theta^m\w) P_\w( X^m_{nt} \leq\fl{nt \vp} + r
The task is to
show that $\frac{1}{\sqrt{n}} W_n(t,r)$ can be approximated by $\frac
{1}{\sqrt{n}} Z_{nt}(\omega)$ uniformly in both $r\in[-R,R]$ and
$t\in
[0,T]$ with probability tending to one. The main work goes
toward approximation uniformly in $t\in[0,T]$ for a fixed $r$.
Uniformity in $r\in[-R,R]$ then comes easily
at the end of this section, completing the proof of Theorem \ref{QMCurrent}.

Before the main work we prove two lemmas that remove a few technical
difficulties. One technical difficulty is presented by small times $t$.
For any fixed $\delta>0$ and $t\geq\delta$ we will use the quenched central
limit theorem to approximate the probabilities in the definition of
$W_n(t,r)$. However, we cannot do this approximation for arbitrarily
small $t$ all at once.
The following lemma will be used later to handle the small values of $t$.
%A difficulty arises in proving process level convergence of $\frac{1}{
%later to overcome this difficulty.
%
\begin{lem}\label{supWnt}
There exists a constant $C<\infty$ such that, for any $r\in\mathbb R$
and $\delta>0$,
\[
\limsup_{n\rightarrow\infty} \frac{1}{\sqrt{n}} E_P \Bigl[ \sup
_{t\in[0,\delta]}
|W_{n}(t,r)| \Bigr] \leq C\sqrt{\delta}.
\]
\end{lem}
\begin{pf}
The triangle inequality implies that
%
%
%e25 ###
\begin{eqnarray}\label{triineq}
&&\frac{1}{\sqrt{n}}E_P \Bigl[ \sup_{t\in[0,\delta]} |W_{n}(t,r) |
\Bigr]
\nonumber
\\[-8pt]
\\[-8pt]
\nonumber
&&\qquad\leq\frac{1}{\sqrt{n}} E_P [ |W_{n}(0,r)| ] + \frac
{1}{\sqrt
{n}} E_P \Bigl[ \sup_{t\in[0,\delta]} |W_{n}(t,r) - W_n(0,r) |
\Bigr].
\end{eqnarray}
For $r>0$,
\[
\frac{1}{\sqrt{n}} W_n(0,r) = \frac{1}{\sqrt{n}} \sum_{0<m\leq
r\sqrt
{n}} E_\omega(\eta_0(m)) - \mu r = \frac{1}{\sqrt{n}} \sum
_{0<m\leq r\sqrt
{n}} \bar{\mu}(\theta^m\omega) - \mu r .
\]
A similar equality holds for $r\leq0$. Therefore, the ergodic theorem
implies that
% \lim_{n\ra\infty} E_P [ |W_{n}(0,r) - \dc r\sqrt{n} | ] =
%0,
the first term on the right-hand side of \eqref{triineq} vanishes as
$n\rightarrow\infty$, and so it remains only to show that
%
%
%e26 ###
\begin{equation}\label{WtrW0r}
\limsup_{n\rightarrow\infty} \frac{1}{\sqrt{n}} E_P \Bigl[ \sup
_{t\in[0,\delta]}
|W_{n}(t,r) - W_n(0,r) | \Bigr] \leq C\sqrt{\delta}.
\end{equation}
Recalling \eqref{oldcurrent} and the fact that $W_n(t,r)= E_\omega Y_n(t,r)
- \mu r \sqrt{n}$, we obtain that
\begin{eqnarray*}
&& W_n(t,r)-W_n(0,r) \\
&&\qquad= \sum_{m> r\sqrt{n}} E_\omega[ \eta_0(m) ] P_\omega\bigl( X^m_{nt}
\leq nt\mathrm{v}_P+
r\sqrt{n} \bigr) \\
&&\qquad\quad {}- \sum_{m\leq r\sqrt{n}}E_\omega[ \eta_0(m) ] P_\omega
\bigl( X^m_{nt} >
nt\mathrm{v}_P+ r\sqrt{n} \bigr).
\end{eqnarray*}
Therefore,
\begin{eqnarray*}
&& \sup_{t\in[0,\delta]} |W_n(t,r)-W_n(0,r)| \\
%& \leq\sum_{m>r\sqrt{n}} E_\w[ \eta_0(m) ] \sup_{t\in[0,\d]} P_\w(
%X_{nt}^m \leq nt\vp+ r\sqrt{n} ) + \sum_{m\leq r\sqrt{n}} E_\w[
&&\qquad\leq\sum_{m>r\sqrt{n}} E_\omega[ \eta_0(m) ] \sup_{t\in
[0,\delta]} P_{\theta
^m \omega}\bigl( X_{nt} - nt\mathrm{v}_P\leq r\sqrt{n} - m\bigr) \\
&&\qquad\quad{} + \sum_{m\leq r\sqrt{n}} E_\omega[ \eta_0(m) ] \sup_{t\in
[0,\delta]}
P_{\theta^m\omega}\bigl(X_{nt} - nt\mathrm{v}_P> r\sqrt{n} - m\bigr) \\
&&\qquad\leq\sum_{m>r\sqrt{n}} E_\omega[ \eta_0(m) ] P_{\theta^m \omega
}\Bigl( \inf
_{t\in[0,\delta]} ( X_{nt} - nt\mathrm{v}_P) \leq
r\sqrt{n} - m \Bigr)
\\
&&\qquad\quad {} + \sum_{m\leq r\sqrt{n}} E_\omega[ \eta_0(m) ] P_{\theta
^m\omega}\Bigl(
\sup_{t\in[0,\delta]} ( X_{nt} - nt\mathrm{v}_P) >
r\sqrt{n} - m \Bigr).
\end{eqnarray*}
Then, the shift invariance of $P$ and Assumption \ref{ICasm} imply that
\begin{eqnarray*}
&& E_P\Bigl\{ \sup_{t\in[0,\delta]} |W_n(t,r)-W_n(0,r)| \Bigr\}\\
&&\qquad\leq E_P\biggl\{ E_\omega[ \eta_0(0) ] \biggl[ \sum
_{m>r\sqrt{n}}
P_{\omega}\Bigl( \inf_{t\in[0,\delta]} ( X_{nt} - nt\mathrm
{v}_P) \leq
r\sqrt{n} - m \Bigr)
\\
&&\qquad\hspace*{71pt}\quad {}+ \sum_{m\leq r\sqrt{n}} P_{\omega}\Bigl( \sup_{t\in[0,\delta]}
( X_{nt} -
nt\mathrm{v}_P) > r\sqrt{n} - m \Bigr) \biggr] \biggr\} \\
&&\qquad \leq E_P\Bigl\{ E_\omega[ \eta_0(0)] \Bigl[ E_\omega
\Bigl(\sup_{t\in
[0,\delta]} ( X_{nt} - nt\mathrm{v}_P)^-\Bigr)\\
&&\qquad\quad \hspace*{68pt}{} +
E_\omega\Bigl(\sup_{t\in
[0,\delta]} ( X_{nt} - nt\mathrm{v}_P)^+ \Bigr) + 1
\Bigr] \Bigr\} \\
% + (r\sqrt{n}-\fl{r\sqrt{n}} ) ] \} \\
&&\qquad \leq2 E_P \Bigl\{ E_\omega[ \eta_0(0) ] E_\omega\Bigl(\sup
_{t\in[0,\delta
]} | X_{nt} - nt\mathrm{v}_P| \Bigr) \Bigr\} + \mu.
\end{eqnarray*}
The Cauchy--Schwarz inequality, along with Assumption \ref{ICMasm} and
Proposition \ref{UIprop}, implies that the right-hand side is bounded
above by $C\sqrt{n\delta} + \mu$.
Dividing by $\sqrt{n}$ and taking $n\rightarrow\infty$, we obtain
\eqref{WtrW0r}.
%By the Cauchy-Schwartz inequality we obtain that
% \E_P \{ E_\w[ \eta_0(0) ] E_\w(\sup_{t\in[0,T]} |
%X_{nt} - nt\vp%| ) \}
%&\leq\E_P \{ E_\w[ \eta_0(0)]^2 \}^{1/2} \E_P \{ %E_\w
%(\sup_{t\in[0,T]} | X_{nt} - nt\vp| )^2
%&\leq C \E[ \sup_{t\in[0,T]} | X_{nt} - nt\vp|^2
%]^{1/2} \leq%C' \sqrt{nT}.
\end{pf}

A second technical difficulty in the analysis of $W_n(t,r)$ is
restricting the sums in the definition of $W_n(t,r)$ to $[-a(n)\sqrt
{n},a(n)\sqrt{n}]$, where $a(n)$ is some sequence tending to $\infty$
slowly (to be specified later, but at least slower than any polynomial
in $n$).
Let $W_n(t,r) = W_{n,1}(t,r) + W_{n,2}(t,r)$, where
\begin{eqnarray*}
W_{n,1}(t,r) &=& \sum_{m=1}^{\lfloor a(n)\sqrt{n} \rfloor}
E_\omega[\eta
_0(m)] P_\omega\bigl(X^m_{nt} \leq nt\mathrm{v}_P+ r\sqrt{n} \bigr) \\
&&{} - \sum_{m=-\lfloor a(n)\sqrt{n} \rfloor+1}^0
  E_\omega[\eta
_0(m)] P_\omega\bigl(X^m_{nt} > nt\mathrm{v}_P+ r\sqrt{n}\bigr) - \mu r \sqrt{n}.
\end{eqnarray*}
%
%W_{n,2}(t,r) &=    \sum_{m=\fl{ a(n)\sqrt{n}}+1}^\infty
%   E_\w[\eta_0(m)] P_\w(X^m_{nt} \leq nt\vp+ r\sqrt{n}) \\
%& - \sum_{m=-\infty}^{-\fl{a(n)\sqrt{n}}}      E_\w[
%and $a(n)$ is some sequence tending to $\infty$ slowly (to be
%specified later, but at least slower than any polynomial in $n$).
The next lemma implies that the main contributions to $W_n(t,r)$ come
from $W_{n,1}(t,r)$.
\begin{lem} \label{Wn2}
For any $\varepsilon>0$, $T<\infty$ and $r\in\mathbb R$,
\[
\lim_{n\rightarrow\infty}P\biggl( \sup_{t\in[0,T]} \frac{1}{\sqrt{n}}|
W_{n,2}(t,r)| \geq\varepsilon\biggr) = 0.
\]
\end{lem}
\begin{pf}
It is enough to show that $E_P |\sup_{t\in[0,T]} W_{n,2}(t,r)| =
o(\sqrt{n})$.
Similarly to the proof of Lemma \ref{supWnt}, we obtain that
\begin{eqnarray*}
&&\sup_{t\in[0,T]} |W_{n,2}(t,r)|\\
%& \leq\sum_{m>\fl{a(n)\sqrt{n}}} E_\w[\eta_0(m)] \sup_{t\in[0,T]} P_
&&\qquad\leq\sum_{m>\lfloor a(n)\sqrt{n} \rfloor} E_\omega[\eta_0(m)]
P_{\theta^m \omega}\Bigl(
\inf_{t\in[0,T]} \bigl( X_{nt} - nt\mathrm{v}_P- r\sqrt{n}
\bigr) \leq- m
\Bigr) \\
&&\qquad\quad {} + \sum_{m\leq-\lfloor a(n)\sqrt{n} \rfloor} E_\omega[\eta
_0(m)] P_{\theta^m\omega
}\Bigl( \sup_{t\in[0,T]} \bigl( X_{nt} - nt\mathrm{v}_P- r\sqrt
{n} \bigr) > -
m \Bigr) ,
\end{eqnarray*}
and the shift invariance of $P$ implies that
\begin{eqnarray*}
&&\hspace*{-4pt} E_P\Bigl\{ \sup_{t\in[0,T]} |W_{n,2}(t,r)| \Bigr\}\\
&&\hspace*{-4pt}\qquad\leq E_P\biggl\{ E_\omega[ \eta_0(0) ] \biggl[ \sum
_{m>\lfloor a(n)\sqrt{n} \rfloor} P_{\omega}\Bigl( \inf_{t\in
[0,T]} \bigl( X_{nt} - nt\mathrm{v}_P- r\sqrt{n}
\bigr) \leq- m \Bigr)
\\
&&\qquad\hspace*{67pt}\quad {}+ \sum_{m\leq-\lfloor a(n)\sqrt{n} \rfloor} P_{\omega}\Bigl( \sup
_{t\in[0,T]} \bigl(
X_{nt} - nt\mathrm{v}_P- r\sqrt{n} \bigr) > - m \Bigr) \biggr]
\biggr\} \\
&&\hspace*{-4pt}\qquad \leq E_P \Bigl\{ E_\omega[ \eta_0(0)] \Bigl[ E_\omega
\Bigl(\sup_{t\in
[0,T]} \bigl( X_{nt} - nt\mathrm{v}_P-r\sqrt{n} + \bigl\lfloor a(n)\sqrt
{n} \bigr\rfloor
\bigr)^-\Bigr)\\
&&\qquad\hspace*{64pt}\quad {} + E_\omega\Bigl(\sup_{t\in[0,T]} \bigl( X_{nt} -
nt\mathrm{v}_P
-r\sqrt{n} - \bigl\lfloor a(n)\sqrt{n} \bigr\rfloor\bigr)^+ \Bigr)
\Bigr] \Bigr\} \\
&&\hspace*{-4pt}\qquad \leq2E_P \Bigl\{ E_\omega[ \eta_0(0)]\Bigl[ E_\omega
\Bigl( \sup_{t\in
[0,T]} \bigl|X_{nt}- nt\mathrm{v}_P- r\sqrt{n}\bigr|
 \\
 &&\qquad\hspace*{120pt}\quad{}
\times  \mathbf{1} \Bigl\{ \sup_{t\in[0,T]} \bigl|X_{nt} - nt\mathrm
{v}_P- r\sqrt{n} \bigr|\\
&&\qquad\hspace*{213pt}\quad{} \geq a(n) \sqrt{n} \Bigr\}  \Bigr) \Bigr]
\Bigr\}.
\end{eqnarray*}
Let $p=2+\varepsilon$ for some $\varepsilon>0$ satisfying Assumption
\ref{ICMasm}, and
let $1/p+1/q = 1$. Note that $p>2$ implies that $q\in(1,2)$. Then,
H\"older's inequality implies that
\begin{eqnarray*}
&&E_P\Bigl\{ \sup_{t\in[0,T]} |W_{n,2}(t,r)| \Bigr\} \\
&&\qquad\leq C E_P\Bigl\{ E_\omega\Bigl( \sup_{t\in[0,T]} \bigl|X_{nt}-
nt\mathrm{v}_P-
r\sqrt{n}\bigr|\\
&&\hspace*{86pt}\qquad{}\times \mathbf{1}{ \Bigl\{ \sup_{t\in[0,T]} \bigl|X_{nt} -
nt\mathrm{v}_P- r\sqrt{n} \bigr| \geq a(n) \sqrt{n} \Bigr\} } \Bigr)
^q \Bigr\}^{1/q},
\end{eqnarray*}
applying the Cauchy--Schwarz inequality to the inner
expectation
%'
\begin{eqnarray*}
&& \leq C E_P\Bigl\{ E_\omega\Bigl( \sup_{t\in[0,T]} \bigl|X_{nt}-
nt\mathrm{v}_P-
r\sqrt{n}\bigr|^2 \Bigr)^{q/2}\\
&&\hspace*{14pt} \qquad{}  \times
P_\omega\Bigl( \sup_{t\in[0,T]} \bigl|X_{nt} - nt\mathrm{v}_P-
r\sqrt{n} \bigr|
\geq a(n) \sqrt{n} \Bigr)^{q/2} \Bigr\}^{1/q}
\end{eqnarray*}
by H\"older's inequality again and because probabilities are
bounded above by $1$
%
%& \leq C \E( \sup_{t\in[0,T]} |X_{nt}- nt\vp- r\sqrt{n}|^2
%)^{1/2} \nn\\
%&  \times E_P \{ P_\w( \sup_{t\in[0,T]}
%|X_{nt} - nt\vp- r\sqrt{n} | \geq a(n) \sqrt{n}
%)^{q/(2-q)} \}^{(2-q)/2q} \nn\\
\begin{eqnarray} \label{EWn2upper}
& \leq& C \mathbb{E}\Bigl( \sup_{t\in[0,T]} \bigl|X_{nt}- nt\mathrm
{v}_P- r\sqrt{n}\bigr|^2
\Bigr)^{1/2}
\nonumber
\\[-8pt]
\\[-8pt]
\nonumber
&&{}\times \mathbb{P}\Bigl( \sup_{t\in[0,T]} \bigl|X_{nt} -
nt\mathrm{v}_P- r\sqrt
{n} \bigr| \geq a(n) \sqrt{n} \Bigr)^{(2-q)/2q}.
\end{eqnarray}
Proposition \ref{UIprop} implies that (for a fixed $T<\infty$ and
$r\in
\mathbb R$) the first term on \eqref{EWn2upper} is $\mathcal{O}(\sqrt
{n})$, and the
averaged functional central limit theorem [part (3) of Theorem~\ref
{QCLTthm}] implies that the last term in \eqref{EWn2upper} vanishes as
$n\rightarrow\infty$. This completes the proof of the lemma.
\end{pf}

%The main result of this section is the following Proposition.
The majority of this section is devoted the the proof of the following
proposition which is a slightly weaker version of Theorem \ref{QMCurrent}.
\begin{prop}\label{Wn1}
For any $\varepsilon>0$, $T<\infty$ and $r\in\mathbb R$,
\[
\lim_{n\rightarrow\infty} P\biggl( \sup_{t\in[0,T]} \frac
{1}{\sqrt{n}} |
W_{n}(t,r) - \mu Z_{nt}(\omega) | \geq\varepsilon\biggr) = 0.
\]
Therefore, $\frac{1}{\sqrt{n}} W_{n}(\cdot,r)$ converges in
distribution to $ \mu\sigma_2 W( \cdot)$, where $W(\cdot)$ is a standard
Brownian motion.
\end{prop}

\begin{pf}
For any $\delta>0$,
\begin{eqnarray}
& &\mathbb{P}\biggl( \sup_{t\in[0,T]} \frac{1}{\sqrt{n}} |
W_{n}(t,r) - \mu
Z_{nt}(\omega) | \geq\varepsilon\biggr) \nonumber\\
&&\qquad \leq\mathbb{P}\biggl( \sup_{t\in[0,\delta]} | W_n(t,r) |
\geq\frac{\varepsilon}{2}
\sqrt{n} \biggr) + \mathbb{P}\biggl( \sup_{t\in[0,\delta]} \mu|
Z_{nt}(\omega) | \geq
\frac{\varepsilon}{2} \sqrt{n} \biggr) \nonumber\\
&&\qquad\quad {}  + \mathbb{P}\biggl( \sup_{t\in[\delta,T]} \frac
{1}{\sqrt{n}} |
W_{n}(t,r) - \mu Z_{nt}(\omega) | \geq\varepsilon\biggr)
\nonumber\\
&&\qquad \leq\frac{2}{\varepsilon\sqrt{n}} \mathbb{E}\Bigl[\sup
_{t\in[0,\delta]} | W_n(t,r)
|\Bigr] + \mathbb{P}\biggl( \sup_{t\in[0,\delta]} \mu|
Z_{nt}(\omega) | \geq\frac{\varepsilon
}{2} \sqrt{n} \biggr) \label{tlessd} \\
&&\qquad\quad {}  + \mathbb{P}\biggl( \sup_{t\in[\delta,T]} \frac
{1}{\sqrt{n}} |
W_{n,2}(t,r) | \geq\frac{\varepsilon}{2} \biggr) \label
{tlessda} \\
&&\qquad\quad {}  + \mathbb{P}\biggl( \sup_{t\in[\delta,T]} \frac
{1}{\sqrt{n}} |
W_{n,1}(t,r) - \mu Z_{nt}(\omega) | \geq\frac{\varepsilon
}{2} \biggr). \label{dtT}
\end{eqnarray}
Letting $n\rightarrow\infty$, Lemma \ref{supWnt} and the fact that
$Z_{nt}(\omega
)/\sqrt{n}$ converges to Brownian motion imply that the two terms in
\eqref{tlessd} can be made arbitrarily small by taking $\delta
\rightarrow0$. Also,
Lemma \ref{Wn2} implies that the term in \eqref{tlessda} vanishes as
$n\rightarrow\infty$. Thus, it is enough to show that, for any
$\delta>0$, \eqref
{dtT} vanishes as $n\rightarrow\infty$.
%To this end, first note that
%P_\w( X_{nt}^m \leq\fl{nt\vp}) &= P_{\theta^m\w}( X_{nt} \leq nt\vp-
%m )\\
%& = P_{\theta^m\w}( \frac{X_{nt} - nt\vp+ Z_{nt}(\theta^m\w)}{
%Theorem \ref{QCLTthm} implies that this last probability is
%approximately
%$\Phi_{\s_1^2 t}( \frac{Z_{nt}(\theta^m\w)-m}{\sqrt{n}} )$
%for $n$ large. Similarly, it should be true that $P_\w( X_{nt}^m >
%The following Lemma makes this approximation precise.
For this, we need the following lemmas whose proofs we defer for now.\vspace*{-1pt}
\begin{lem} \label{WtildeW}
Let
\begin{eqnarray*}
\widetilde{W}_{n,1}(t,r) &:=& \sum_{m=1}^{\lfloor a(n)\sqrt{n}
\rfloor} E_\omega( \eta_0(m))
\Phi_{\sigma_1^2 t}\biggl( \frac{Z_{nt}(\theta^m\omega)-m}{\sqrt
{n}} + r \biggr)
\\[-1pt]
&&{} - \sum_{m=-\lfloor a(n)\sqrt{n} \rfloor+1}^0     E_\omega
( \eta_0(m)) \Phi
_{\sigma_1^2 t}\biggl( - \frac{Z_{nt}(\theta^m\omega)-m}{\sqrt{n}}
- r \biggr)\\
&&{} -
\mu r \sqrt{n}.
\end{eqnarray*}
Then, for any $\varepsilon>0$, $r\in\mathbb R$ and $0<\delta
<T<\infty$,
\[
\lim_{n\rightarrow\infty} P\biggl( \sup_{t \in[\delta,T]} \frac
{1}{\sqrt{n}}
|W_{n,1}(t,r) - \widetilde{W}_{n,1}(t,r)| \geq\varepsilon
\biggr) = 0.
\]
\end{lem}
%
%The next step in the proof of Proposition \ref{Wn1} is to replace
%$Z_{nt}(\theta^m\w)$ by $Z_{nt}(\w)$.
% for all $|m|\leq a(n)\sqrt{n}$.
%
\begin{lem} \label{tWhW}
Let
\begin{eqnarray*}
\widehat{W}_{n,1}(t,r) &:=& \sum_{m=1}^{\lfloor a(n)\sqrt{n} \rfloor
} E_\omega( \eta_0(m))
\Phi_{\sigma_1^2 t}\biggl( \frac{Z_{nt}(\omega)-m}{\sqrt{n}} + r
\biggr) \\[-1pt]
&&{} -     \sum_{m=-\lfloor a(n)\sqrt{n} \rfloor+1}^0
E_\omega( \eta
_0(m)) \Phi_{\sigma_1^2 t}\biggl( - \frac{Z_{nt}(\omega)-m}{\sqrt
{n}} - r
\biggr) - \mu r \sqrt{n}.
\end{eqnarray*}
Then, for any $\varepsilon>0$, $r\in\mathbb R$ and $0<\delta
<T<\infty$,
\[
\lim_{n\rightarrow\infty} P\biggl( \sup_{t\in[\delta,T]} \frac
{1}{\sqrt{n}}|\widetilde
{W}_{n,1}(t,r) - \widehat{W}_{n,1}(t,r)| \geq\varepsilon
\biggr) = 0.
\]
\end{lem}
\begin{lem} \label{hWbW}
Let
\begin{eqnarray*}
\overline{W}_{n,1}(t,r) &:=& \sum_{m=1}^{\lfloor a(n)\sqrt{n}
\rfloor} \mu\Phi_{\sigma_1^2
t}\biggl( \frac{Z_{nt}(\omega)-m}{\sqrt{n}} + r \biggr) \\
&&{} -     \sum_{m=-\lfloor a(n)\sqrt{n} \rfloor+1}^0
\mu\Phi_{\sigma
_1^2 t}\biggl( - \frac{Z_{nt}(\omega)-m}{\sqrt{n}} - r \biggr) -
\mu r \sqrt{n}.
\end{eqnarray*}
Then, for any $\varepsilon>0$, $r\in\mathbb R$ and $0<\delta
<T<\infty$,
\[
\lim_{n\rightarrow\infty} P\biggl( \sup_{t\in[\delta,T]} \frac
{1}{\sqrt{n}}|\widehat
{W}_{n,1}(t) - \overline{W}_{n,1}(t)| \geq\varepsilon
\biggr) = 0.
\]
\end{lem}

Assuming for now Lemmas \ref{WtildeW}, \ref{tWhW} and \ref{hWbW}, to
finish the proof of Proposition~\ref{Wn1}, it remains to compare
$\overline
{W}_{n,1}(t,r)$ with $\mu Z_{nt}(\omega)$.
%For this, we introduce parameters $R$ and $L$ which will eventually be
%made arbitrarily large after first letting $n\ra\infty$.
%Then,
% \frac{1}{\sqrt{n}} \wh{W}_{n,1}(t) &= \frac{1}{\sqrt{n}} \sum_{
% & - \frac{1}{\sqrt{n}} \sum_{\ell=-\fl{R/L} }^{-1} \sum_{m=\fl{
%& + \frac{1}{\sqrt{n}}\sum_{m=R\sqrt{n}}^{\fl{a(n)\sqrt{n}}} E_
%) \\
%& - \frac{1}{\sqrt{n}}\sum_{m=-\fl{a(n)\sqrt{n}}}^{-R\sqrt{n}} E_
%)
Since $\Phi_{\sigma_1^2t}(\cdot)$ is strictly increasing and bounded above
by 1, we have using a Riemann sum approximation that, for any $t\in[0,T]$,
%
%
%e27 ###
\begin{eqnarray} \label{RSapprox}
\hspace*{18pt}&&\biggl| \frac{\overline{W}_{n,1}(t,r)}{\sqrt{n}} + \mu r
\nonumber\\
&&{} - \mu\int
_0^{a(n)}
\biggl(\Phi_{\sigma_1^2 t}\biggl( \frac{Z_{nt}(\omega)}{\sqrt{n}} + r -x
\biggr) - \Phi
_{\sigma_1^2 t}\biggl( - \frac{Z_{nt}(\omega)}{\sqrt{n}} - r -x
\biggr) \biggr)\,dx
\biggr| \\
&&\qquad\leq\frac{2\mu}{\sqrt{n}}.\nonumber
\end{eqnarray}
%
%Therefore, we are reduced to showing that $\int_0^{a(n)} \Phi_{\s_1^2
%t}( \frac{Z_{nt}(\w)}{\sqrt{n}} -x ) - \Phi_{\s_1^2 t}
%( - \frac{Z_{nt}(\w)}{\sqrt{n}} -x ) dx $ is approximately
%equal to $\frac{1}{\sqrt{n}} Z_{nt}(\w)$.
It is an easy exercise in calculus to show that, for any $z\in\mathbb
R$ and $A>0$,
\[
\int_0^{A} \bigl(\Phi_{\alpha^2}( z -x ) - \Phi_{\alpha
^2}( -z -x
)\bigr)\, dx = z + \Psi_{\alpha^2}(A+z) - \Psi_{\alpha^2}(A-z),
\]
where $\Psi_{\alpha^2}(x)$ is defined in \eqref{Psidef}.
%$\Psi_{\a^2}(x) := \a^2 \phi_{\a^2}(x)-x\Phi_{\a^2}(-x)$ and $\phi_{
%distribution with variance $\a^2$.
Therefore,
\begin{eqnarray*}
&&\int_0^{a(n)} \biggl(\Phi_{\sigma_1^2 t}\biggl( \frac{Z_{nt}(\omega
)}{\sqrt{n}} +r -x
\biggr) - \Phi_{\sigma_1^2 t}\biggl( - \frac{Z_{nt}(\omega)}{\sqrt
{n}} - r -x
\biggr)\biggr)\, dx \\
&&\qquad = \frac{Z_{nt}(\omega)}{\sqrt{n}} + r + \Psi_{\sigma_1^2
t}
\biggl(a(n)+\frac{Z_{nt}(\omega)}{\sqrt{n}} + r \biggr)\\
&&\qquad\quad {} - \Psi_{\sigma
_1^2 t}
\biggl(a(n)-\frac{Z_{nt}(\omega)}{\sqrt{n}} - r \biggr).
\end{eqnarray*}
Recalling \eqref{RSapprox}, this implies that for $\varepsilon>0$ and $n$
sufficiently large,
\begin{eqnarray}\label{Psireduce}
&& P\biggl( \sup_{t\in[0,T]} \biggl| \frac{\overline
{W}_{n,1}(t,r)}{\sqrt{n}} -
\frac{Z_{nt}(\omega)}{\sqrt{n}} \biggr| \geq\varepsilon\biggr)
\nonumber\\
&&\qquad \leq P\biggl( \sup_{t\in[0,T]} \biggl| \Psi_{\sigma_1^2
t}
\biggl(a(n)+\frac{Z_{nt}(\omega)}{\sqrt{n}}+ r \biggr)\\
&&\qquad\hspace*{42pt}\quad{} - \Psi_{\sigma
_1^2 t}
\biggl(a(n)-\frac{Z_{nt}(\omega)}{\sqrt{n}} - r \biggr) \biggr| \geq
\frac{\varepsilon}{2
\mu} \biggr).\nonumber
%& \leq P( \sup_{t\in[0,T]} | \Psi_{\s_1^2 t}(a(n)+
%& \leq P( \sup_{t\in[0,T]} | \frac{Z_{nt}(\w)}{
\end{eqnarray}
A simple calculation shows that $\Psi_{\alpha^2}'(x) = - \Phi
_{\alpha^2}(-x) <
0$, and so $\Psi_{\alpha^2}(x)$ is decreasing in $x$.
Another direct calculation shows that $\frac{d}{d\alpha} \Psi
_{\alpha^2} (x) =
\alpha\phi_{\alpha^2}(x) >0$. Thus, $\Psi_{\alpha^2}(x)$ is
increasing in $\alpha$.
Thus, if $|Z_{nt}(\omega)| \leq a(n)\sqrt{n}/2$ and $t\leq T$,
\begin{eqnarray*}
&&\biggl| \Psi_{\sigma_1^2 t}\biggl(a(n)+\frac{Z_{nt}(\omega)}{\sqrt
{n}} + r
\biggr) - \Psi_{\sigma_1^2 t}\biggl(a(n)-\frac{Z_{nt}(\omega)}{\sqrt{n}}
- r \biggr)
\biggr|\\
 &&\qquad\leq2 \Psi_{\sigma_1^2 t}\bigl(a(n)/2 - |r| \bigr) \\
&&\qquad\leq2 \Psi_{\sigma_1^2 T}\bigl( a(n)/2 - |r| \bigr).
\end{eqnarray*}
%
%Another direct calculation shows that $\frac{d}{d\a} \Psi_{\a^2} (x) =
%decreases as $\a$ increases, and so
% \sup_{t\in[0,T]} \Psi_{\s_1^2 t}(a(n)- | \frac{Z_{nt}(\w)}{
%Since $n^{-1/2} Z_{nt}(\w)$ converges in distribution to $W_{\s_2^2
%t}$, we need only show that the last two terms converge to zero in
%probability.
Since $\lim_{x\rightarrow\infty} \Psi_{\sigma_1^2 T}(x) = 0$,\vspace*{1pt} then
$\Psi_{\sigma_1^2
T}(a(n)/2-|r|) < \frac{\varepsilon}{2}$ for all $n$ large enough. Thus,
recalling \eqref{Psireduce}, we obtain that, for any $\varepsilon>0$
and $n$
sufficiently large,
\[
P\biggl( \sup_{t\in[0,T]} \biggl| \frac{\overline
{W}_{n,1}(t,r)}{\sqrt{n}} -
\frac{Z_{nt}(\omega)}{\sqrt{n}} \biggr| \geq\varepsilon\biggr)
%P( \sup_{t\in[0,T]} | \Psi_{\s_1^2 t}(a(n)+\frac{Z_{nt}(
\leq P\biggl( \sup_{t\in[0,T]}\biggl | \frac{Z_{nt}(\omega)}{\sqrt
{n}}
\biggr| \geq\frac{a(n)}{2} \biggr).
\]
Since $t\mapsto\frac{Z_{nt}(\omega)}{\sqrt{n}}$ converges in distribution
to a Brownian motion, this last probability tends to zero as
$n\rightarrow\infty
$. This completes the proof of Proposition \ref{Wn1}.
\end{pf}

We now return to the proofs of Lemmas \ref{WtildeW}--\ref{hWbW}.
\begin{pf*}{Proof of Lemma \protect\ref{WtildeW}}
Let
\begin{eqnarray*}
D(n,\omega) &:=& \sup_{x\in\mathbb R} \biggl| P_\omega\biggl( \frac
{X_n - n\mathrm{v}_P+ Z_n(\omega
)}{\sqrt{n}} \leq x \biggr) - \Phi_{\sigma_1^2}(x) \biggr|\quad
\mbox
{and}\\
\bar{D}(n,\omega) &:=& \sup_{k\geq n} D(k,\omega).
\end{eqnarray*}
Theorem \ref{QCLTthm} implies that $\lim_{n\rightarrow\infty} \bar
{D}(n,\omega) =
0$, $P\mathrm{\mbox{-}a.s.}$, and so by the bounded convergence theorem, $\lim
_{n\rightarrow
\infty} E_P[ \bar{D}(n,\omega)^p] = 0$ for any $p>0$. Thus, it is possible
to choose the sequence $a(n)$ tending to infinity slowly enough so that
\[
\lim_{n\rightarrow\infty} a(n) (E_P[ \bar{D}(\delta
n,\omega)^2 ]
)^{1/2} = 0\qquad  \forall\delta>0
\]
[e.g., let $a(n) = (E_P[ \bar{D}(\sqrt{n},\omega)^2
])^{-1/4}$]. The definition of $D(n,\omega)$ implies that, for any $t>0$,
\[
| W_{n,1}(t,r) - \widetilde{W}_{n,1}(t,r) | \leq%
\sum_{m=-\lfloor a(n)\sqrt{n} \rfloor+1}^{\lfloor a(n)\sqrt{n}
\rfloor} E_{\theta^m\omega}( \eta
_0(0)) D(nt,\theta^m\omega) .
\]
Therefore,
\begin{eqnarray*}
&&P\Bigl( \sup_{t\in[\delta,T]} | W_{n,1}(t,r) - \widetilde
{W}_{n,1}(t,r)
| \geq\varepsilon\sqrt{n} \Bigr) \\
&&\qquad\leq P\Biggl( \sup_{t\in[\delta,T]} \sum_{m=-\lfloor a(n)\sqrt{n}
\rfloor+1}^{\lfloor a(n)\sqrt{n} \rfloor} E_{\theta^m\omega}( \eta
_0(0)) D(nt,\theta^m\omega) \geq\varepsilon\sqrt
{n} \Biggr) \\
&&\qquad\leq P\Biggl( \sum_{m=-\lfloor a(n)\sqrt{n} \rfloor+1}^{\lfloor
a(n)\sqrt{n} \rfloor}
E_{\theta^m\omega}( \eta_0(0)) \bar{D}(\delta n,\theta^m\omega)
\geq\varepsilon\sqrt{n}
\Biggr) \\
%&\leq\frac{1}{\e\sqrt{n}} \sum_{m=-\fl{a(n)\sqrt{n}}+1}^{\fl{a(n)
&&\qquad\leq\frac{2 a(n)}{\varepsilon} E_P[E_{\omega}( \eta_0(0)) \bar
{D}(\delta n,\omega)] \\
&&\qquad\leq\frac{2 a(n)}{\varepsilon} ( E_P[(E_\omega\eta_0(0))^2]
)^{1/2}
( E_P[ \bar{D}(\delta n,\omega)^2] )^{1/2},
%&=\frac{2 a(n)}{\e} \| E_\w(\eta_0(0))] \|_{L^2(P)} \|
\end{eqnarray*}
where the next to last inequality follows from Chebyshev's inequality
and the shift invariance of $P$.
%An application of the Cauchy-Schwartz inequality, along with
%Assumption \ref{ICMasm} and the conditions on the rate of growth of
%$a(n)$ implies that this last term vanishes as $n\ra\infty$.
Our choice of the sequence $a(n)$ ensures that this last term vanishes
as $n\rightarrow\infty$.
\end{pf*}

\begin{pf*}{Proof of Lemma \protect\ref{tWhW}}
Note that the mean value theorem implies %since $\sup_{x} \Phi_{
\begin{eqnarray*}
| \Phi_{\sigma_1^2t}(x) - \Phi_{\sigma_1^2 t}(y) |
&\leq&\Bigl( \sup
_{z\in\mathbb R} \Phi_{\sigma_1^2t}'(z) \Bigr) |x-y|\\
& =& \frac
{1}{\sigma_1\sqrt{2\pi
t}}|x-y|\qquad  \forall x,y \in\mathbb R.
\end{eqnarray*}
Therefore,
\begin{eqnarray*}
&&\sup_{t\in[\delta,T]} |\widetilde{W}_{n,1}(t,r) - \widehat
{W}_{n,1}(t,r)| \\
&&\qquad\leq\sup
_{t\in[\delta,T]}
\sum_{m=-\lfloor a(n)\sqrt{n} \rfloor+1}^{\lfloor a(n)\sqrt{n}
\rfloor} \bar\mu(\theta^m \omega)
\frac{1}{\sigma_1\sqrt{2\pi t}} \biggl| \frac{Z_{nt}(\theta
^m\omega) - Z_{nt}(\omega
)}{\sqrt{n}} \biggr| \\
&&\qquad\leq\frac{2 a(n)}{\sigma_1 \sqrt{2\pi\delta}} \sup_{t\in
[\delta,T]} \max_{| m|
\leq a(n)\sqrt{n}} | Z_{nt}(\theta^m\omega) - Z_{nt}(\omega)
| \\
&&\qquad\quad {} \times\Biggl( \frac{1}{2 a(n)\sqrt{n}} \sum_{m=-\lfloor
a(n)\sqrt{n} \rfloor+1}^{\lfloor a(n)\sqrt{n} \rfloor} \bar\mu
(\theta^m \omega) \Biggr).
\end{eqnarray*}
The ergodic theorem implies that the averaged sum on the last line
converges to $\mu$, $P\mathrm{\mbox{-}a.s.}$ Thus, to finish the proof of the lemma, it
is enough to show that
\[
\lim_{n\rightarrow\infty} P\biggl( \sup_{t\in[\delta,T]} \max_{|
m| \leq a(n)\sqrt
{n}} | Z_{nt}(\theta^m\omega) - Z_{nt}(\omega) | \geq
\frac{\varepsilon\sqrt
{n}}{a(n)} \biggr) = 0\qquad   \forall\varepsilon>0.
\]
Since $Z_{nt}(\theta^m\omega) = h(m+\lfloor nt\mathrm{v}_P \rfloor
,\omega)-h(m,\omega)$,
\[
|Z_{nt}(\theta^m\omega) - Z_{nt}(\omega)| \leq|h(m,\omega)| +
|h(m+\lfloor nt\mathrm{v}_P \rfloor,\omega)
- h(\lfloor nt\mathrm{v}_P \rfloor,\omega) |.
\]
Thus,
\begin{eqnarray*}
&& \sup_{t\in[\delta,T]} \max_{| m| \leq a(n)\sqrt{n}} |
Z_{nt}(\theta
^m\omega) - Z_{nt}(\omega) | \\
&&\qquad  \leq2 \max_{x\in[0,nT]} \max_{1\leq m \leq a(n)\sqrt{n} }
|h(x+m,\omega) - h(x,\omega) | \\
%& \leq4 \max_{0 \leq i\leq\sqrt{n}T/a(n)} \max_{1 \leq l,m \leq a(n)
&&\qquad  \leq6 \max_{0 \leq i\leq\sqrt{n}T/a(n)} \max_{1 \leq m
\leq
a(n)\sqrt{n}} \bigl|h\bigl(i \bigl\lfloor a(n)\sqrt{n} \bigr\rfloor+ m,\omega\bigr) - h\bigl(i
\bigl\lfloor a(n)\sqrt{n} \bigr\rfloor,\omega\bigr)\bigr|.
\end{eqnarray*}
This implies that
\begin{eqnarray*}
&& P\biggl( \sup_{t\in[\delta,T]} \max_{| m| \leq a(n)\sqrt{n}}
|
Z_{nt}(\theta^m\omega) - Z_{nt}(\omega) | \geq\frac
{\varepsilon\sqrt{n}}{a(n)}
\biggr)\\
&&\qquad \leq P\biggl( \max_{0 \leq i\leq\sqrt{n}T/a(n)} \max_{1
\leq m
\leq a(n)\sqrt{n}} \bigl|h\bigl(i \bigl\lfloor a(n)\sqrt{n} \bigr\rfloor+ m,\omega\bigr) -
h\bigl(i \bigl\lfloor a(n)\sqrt{n} \bigr\rfloor,\omega\bigr)\bigr| \\
&&\hspace*{288pt}\qquad\geq\frac{\varepsilon
\sqrt{n}}{6a(n)} \biggr)\\
&&\qquad\leq\frac{\sqrt{n}T}{a(n)} P\biggl( \max_{1 \leq m \leq
a(n)\sqrt{n}} |h(m,\omega)| \geq\frac{\varepsilon\sqrt{n}}{6a(n)}
\biggr),
\end{eqnarray*}
where the last inequality is from a union bound and the shift
invariance of $P$. %Goldsheid has shown \cite[Lemma 4]{gQCLT} that
%there exists an $\eta> 0$ and a constant $C<\infty$ such that
% E_P [ \sup_{1\leq k \leq n} |h(k,\w)|^{2+2\eta} ] \leq C
%n^{1+\eta},  \forall n\in\N.
Recalling Lemma \ref{Goldsheid}, there exist constants $C,\eta>0$ such
that, for any fixed $\varepsilon>0$ and $0<\delta<T<\infty$,
%Applying Chebyshev's inequality to \eqref{tWsuph} then implies that
%for any fixed $\e>0$ and $0<\d<T<\infty$,
%
\begin{eqnarray*}
&& P\biggl( \sup_{t\in[\delta,T]} \max_{| m| \leq a(n)\sqrt{n}} |
Z_{nt}(\theta^m\omega) - Z_{nt}(\omega) | \geq\frac
{\varepsilon\sqrt{n}}{a(n)}
\biggr)\\
&&\qquad\leq\frac{\sqrt{n}T}{a(n)} \biggl( \frac{6a(n)}{\varepsilon\sqrt
{n}}
\biggr)^{2+2\eta} C \bigl( a(n)\sqrt{n} \bigr)^{1+\eta} \\
&&\qquad= \mathcal{O}( a(n)^{2+3\eta} n^{-\eta/2} ).
\end{eqnarray*}
Since $a(n)$ grows slower than polynomially in $n$, this last term
vanishes as $n\rightarrow\infty$.
%\rightqed
\end{pf*}

\begin{pf*}{Proof of Lemma \protect\ref{hWbW}}
For any integer $R$ let
\begin{eqnarray*}
\widehat{W}_{n,1}^R(t,r) &:=& \sum_{m=1}^{\lfloor R \sqrt{n} \rfloor
} E_\omega( \eta_0(m))
\Phi_{\sigma_1^2 t}\biggl( \frac{Z_{nt}(\omega)-m}{\sqrt{n}} + r
\biggr) \\
&&{} -     \sum_{m=-\lfloor R\sqrt{n} \rfloor+1}^0
E_\omega( \eta_0(m))
\Phi_{\sigma_1^2 t}\biggl( - \frac{Z_{nt}(\omega)-m}{\sqrt{n}} - r
\biggr) - \mu
r \sqrt{n}
\end{eqnarray*}
and
\begin{eqnarray*}
\overline{W}{}^R_{n,1}(t,r) &:=& \sum_{m=1}^{ \lfloor R \sqrt{n}
\rfloor} \mu\Phi_{\sigma_1^2
t}\biggl( \frac{Z_{nt}(\omega)-m}{\sqrt{n}} + r \biggr) \\
&&{}-
\sum_{m=-
\lfloor R \sqrt{n} \rfloor+1}^0     \mu\Phi_{\sigma_1^2 t}
\biggl( - \frac{Z_{nt}(\omega
)-m}{\sqrt{n}} - r \biggr) - \mu r \sqrt{n}.
\end{eqnarray*}
Then, it is enough to show that
%
%
%e28 ###
\begin{equation}\label{hWbWR}
\lim_{n\rightarrow\infty} \frac{1}{\sqrt{n}} E_P \Bigl[ \sup
_{t\in[\delta,T]}
|\widehat{W}_{n,1}^R(t,r) - \overline{W}{}^R_{n,1}(t,r)|
\Bigr] = 0\qquad
\forall R<\infty,
\end{equation}
and that
%
%
%e29 ###
\begin{eqnarray}\label{hWR}
&& \lim_{R\rightarrow\infty} \limsup_{n\rightarrow\infty} P\biggl(
\sup_{t\in[\delta,T]} \frac
{1}{\sqrt{n}} | \widehat{W}_{n,1}(t,r) - \widehat
{W}_{n,1}^R(t,r) |
\geq\varepsilon\biggr)
\nonumber
\\[-8pt]
\\[-8pt]
\nonumber
&&\qquad = 0\qquad  \forall\varepsilon>0
\end{eqnarray}
and
%
%
%e30 ###
\begin{eqnarray}\label{bWR}
&&\lim_{R\rightarrow\infty} \limsup_{n\rightarrow\infty} P\biggl(
\sup_{t\in[\delta,T]} \frac
{1}{\sqrt{n}} | \overline{W}_{n,1}(t,r) - \overline
{W}{}^R_{n,1}(t,r) |
\geq\varepsilon\biggr)
\nonumber
\\[-8pt]
\\[-8pt]
\nonumber
&&\qquad= 0\qquad  \forall\varepsilon>0.
\end{eqnarray}
%
% P( \sup_{t\in[\d,T]} \frac{1}{\sqrt{n}} | \wh{W}_{n,1}(t) -
%& \leq P( \sup_{t\in[\d,T]} \frac{1}{\sqrt{n}} |
%& + P( \sup_{t\in[\d,T]} \frac{1}{\sqrt{n}} |
%& + P( \sup_{t\in[\d,T]} \frac{1}{\sqrt{n}} |
%). \label{hWbWR}
%It is enough to show that \eqref{hWbWR} vanishes as $n\ra\infty$ for
%any fixed $R$ and that the $\limsup$ as $n\ra\infty$ of \eqref{hWR}
%and \eqref{bWR} can be made arbitrarily small by choosing $R$ large
%enough.

To bound \eqref{hWbWR}, we fix another parameter $L$ and then divide
the interval $(-\lfloor R\sqrt{n} \rfloor, \lfloor R\sqrt{n} \rfloor
]$ into $2 RL$ intervals,
each of length approximately $ \sqrt{n}/L $. For ease of notation, let
$B_{n,L}(\ell) := \{ m \in\mathbb Z\dvtx  \frac{(\ell-1)\sqrt{n}}{L} <
m \leq\frac
{\ell\sqrt{n}}{L} \}$.
%Then $(-\fl{R\sqrt{n}}, \fl{R\sqrt{n}}] \cap\Z= \bigcup_{\ell=
%-RL+1}^{RL} B_{n,L}(\ell).$
Now, for any $m\in B_{n,L}(\ell)$ and $t \in[\delta,T]$,
\begin{eqnarray*}
\Phi_{\sigma_1^2 t} \biggl( \frac{Z_{nt}(\omega) - m}{\sqrt{n}} +
r \biggr) -
\Phi_{\sigma_1^2 t}\biggl ( \frac{Z_{nt}(\omega)}{\sqrt{n}} -
\frac{\ell}{L} + r
\biggr) &\leq&\frac{1}{\sqrt{2\pi t}} \biggl| \frac{m}{\sqrt{n}} -
\frac
{\ell}{L} \biggr|\\
& \leq&\frac{C}{L},
\end{eqnarray*}
where the constant $C$ depends only on $\delta>0$. Thus,\vspace*{-1pt}
\begin{eqnarray*}
\frac{1}{\sqrt{n}} \widehat{W}_{n,1}^R(t,r) &=& \sum_{\ell=
1}^{RL}\biggl(
\frac{1}{\sqrt{n}} \sum_{m\in B_{n,L}(\ell)} \bar\mu(\theta
^m\omega) \biggr)
\Phi_{\sigma_1^2 t} \biggl( \frac{Z_{nt}(\omega)}{\sqrt{n}} -
\frac{\ell}{L} + r
\biggr) \\[-1pt]
&&{} - \sum_{\ell= -RL+1}^{0}\biggl( \frac{1}{\sqrt{n}} \sum
_{m\in
B_{n,L}(\ell)} \bar\mu(\theta^m\omega) \biggr) \Phi_{\sigma_1^2
t}\biggl ( - \frac
{Z_{nt}(\omega)}{\sqrt{n}} + \frac{\ell}{L} - r \biggr)\\[-1pt]
&&{} + \frac{1}{\sqrt{n}} \sum_{m=1}^{\lfloor R \sqrt{n} \rfloor
} \bar\mu(\theta
^m\omega) \mathcal{O}(L^{-1}) \\[-1pt]
&&{}- \frac{1}{\sqrt{n}} \sum
_{m=-\lfloor R\sqrt{n} \rfloor+
1}^{0} \bar\mu(\theta^m\omega) \mathcal{O}(L^{-1}).\vspace*{-1pt}
\end{eqnarray*}
A similar equality also holds for $\overline{W}{}^R_{n,1}(t,r)$ with
$\bar\mu
(\theta^m\omega)$ replaced by $\mu$. Therefore, using the fact that
$\Phi
_{\alpha^2}$ is bounded by 1, we obtain that
\begin{eqnarray*}
&& \sup_{t\in[\delta,T]} \frac{1}{\sqrt{n}} |\widehat
{W}_{n,1}^R(t,r) - \overline
{W}{}^R_{n,1}(t,r)|\\[-1pt]
&&\qquad \leq\sum_{\ell= -RL+1}^{RL}\biggl| \frac{1}{\sqrt{n}} \sum
_{m\in B_{n,L}(\ell)} \bigl( \bar\mu(\theta^m\omega) - \mu\bigr)
\biggr|\\[-1pt]
&&\qquad\quad {}+ \mathcal{O}(L^{-1}) \Biggl( \frac{1}{\sqrt{n}} \sum_{m=-\lfloor
R\sqrt{n} \rfloor+
1}^{\lfloor R \sqrt{n} \rfloor} \bar\mu(\theta^m\omega) + 2 R \mu
\Biggr).
\end{eqnarray*}
Note that we were able to include the supremum over $t$ in the above
inequality since the constant in the $\mathcal{O}(L^{-1})$ term is
valid for
any $t\geq\delta$. Taking expectations of the above with respect to the
measure $P$ and letting $n\rightarrow\infty$, the ergodic theorem
implies that
the first term vanishes and the second term has $\limsup$ less than
$4R\mu\mathcal{O}(L^{-1})$. Thus, taking $L\rightarrow\infty$
proves \eqref{hWbWR}.

To bound \eqref{hWR}, let\vspace*{-1pt}
\[
G_{n,R}:= \biggl\{ \omega\dvtx  \sup_{t\in[\delta,T]} \biggl| \frac
{Z_{nt}(\omega)}{\sqrt
{n}} + r \biggr| \leq\frac{R}{2} \biggr\}.
\]
Since $t\mapsto Z_{nt}(\omega)/\sqrt{n}$ converges to Brownian
motion,
\[
\lim
_{R\rightarrow\infty} \lim_{n\rightarrow\infty} P(G_{n,R}) = 1
\]
for any fixed $r\in\mathbb R$. Thus,
\begin{eqnarray}\label{hWRG}
&& \lim_{R\rightarrow\infty} \limsup_{n\rightarrow\infty} P\biggl(
\sup_{t\in[\delta,T]} \frac
{1}{\sqrt{n}} | \widehat{W}_{n,1}(t,r) - \widehat
{W}_{n,1}^R(t,r) |
\geq\varepsilon\biggr)
\nonumber
\\[-8pt]
\\[-8pt]
\nonumber
&&\qquad \leq\lim_{R\rightarrow\infty}\limsup_{n\rightarrow\infty
} \frac{1}{\varepsilon\sqrt
{n}} E_P\Bigl[\sup_{t\in[\delta,T]} |\widehat{W}_{n,1}(t,r)
- \widehat
{W}_{n,1}^R(t,r) | \mathbf{1}_{G_{n,R}} \Bigr].
\end{eqnarray}
If $\omega\in G_{n,R}$, $|m| \geq R\sqrt{n}$ and $t\leq T$, then
$\Phi_{\sigma
_1^2 t} ( | \frac{Z_{nt}(\omega)}{\sqrt{n}} + r |
- \frac
{|m|}{\sqrt{n}} ) \leq\Phi_{\sigma_1^2 T} ( \frac
{R}{2} - \frac
{|m|}{\sqrt{n}} )$.
Therefore,
\begin{eqnarray*}
&& \frac{1}{\sqrt{n}} E_P\Bigl[\sup_{t\in[\delta,T]}
|\widehat{W}_{n,1}(t) -
\widehat{W}_{n,1}^R(t) | \mathbf{1}_{G_{n,R}} \Bigr] \\
&&\qquad \leq\frac{1}{\sqrt{n}} \sum_{m=\lfloor R\sqrt{n} \rfloor
+1}^{\lfloor a(n)\sqrt{n} \rfloor} \mu\Phi_{\sigma_1^2 T}\biggl(
\frac{R}{2} - \frac{m}{\sqrt{n}}
\biggr)\\
&&\qquad\quad {} + \frac{1}{\sqrt{n}} \sum_{m=-\lfloor a(n)\sqrt{n} \rfloor
+1}^{-\lfloor R\sqrt{n} \rfloor}
\mu\Phi_{\sigma_1^2 T}\biggl( \frac{R}{2} + \frac{m}{\sqrt{n}}
\biggr) \\
&&\qquad \leq\mu\int_R^\infty\Phi_{\sigma_1^2 T}( R/2 - x
) \,dx
+ \mu\int_{-\infty}^{-R} \Phi_{\sigma_1^2 T} ( R/2 + x
) \,dx +
\frac{\mu}{\sqrt{n}},
\end{eqnarray*}
where the last inequality is from a Riemann sum approximation. Since
the integrals in the last line can be made arbitrarily small by taking
$R\rightarrow\infty$, recalling \eqref{hWRG} finishes the proof of
\eqref{hWR}.
The proof of \eqref{bWR} is similar.
\end{pf*}

We conclude this section with the proof of Theorem \ref{QMCurrent}.
\begin{pf*}{Proof of Theorem \protect\ref{QMCurrent}}
To prove Theorem \ref{QMCurrent} from Proposition \ref{Wn1}, we need to
justify the ability to include a supremum over $r\in[-R,R]$ inside the
probability in the statement of Proposition \ref{Wn1}. A simple union
bound implies that we may include a supremum over a finite set of $r$
values inside the probability in the statement of Proposition \ref{Wn1}.
That is, for $N<\infty$ and $r_1,r_2,\ldots, r_N \in\mathbb
R$,\looseness=1
%
%
%e31 ###
\begin{equation}\label{finiter}
\qquad\quad \lim_{n\rightarrow\infty} P\biggl( \max_{k\leq N} \sup_{t\in
[0,T]} \frac{1}{\sqrt
{n}} \bigl| E_\omega Y_{n}(t,r_k) - \mu r_k \sqrt{n} - \mu
Z_{nt}(\omega)
\bigr| \geq\varepsilon\biggr) = 0.
\end{equation}
Now, the definition of $Y_n(t,r)$ implies that $Y_n(t,r)$ is
nondecreasing in $r$.
%We will take advantage of this fact to improve \eqref{finiter} to
%include a supremum over $r\in[-R,R]$.
Therefore, for any fixed $t$, $E_\omega Y_n(t,r) - \mu Z_{nt}(\omega)$ is
nondecreasing in $r$.
Choose $-R=r_1 < r_2 < \cdots< r_{N-1} < r_N = R$ such that
$r_{k+1}-r_k \leq\frac{\varepsilon}{2\mu}$ for $k=1,\ldots, N-1$.
Then, if $r\in[r_k,r_{k+1}]$,
\begin{eqnarray*}
&& \bigl\{ \bigl| E_\omega Y_{n}(t,r) - \mu r \sqrt{n} - \mu
Z_{nt}(\omega)
\bigr| \geq\varepsilon\sqrt{n} \bigr\} \\
&&\qquad\subset\biggl\{ \bigl| E_\omega Y_{n}(t,r_k) - \mu r_k \sqrt{n} -
\mu
Z_{nt}(\omega) \bigr| \geq\frac{\varepsilon}{2}\sqrt{n}
\biggr\} \\
&&\hspace*{12pt}\qquad{}  \cup\biggl\{ \bigl| E_\omega Y_{n}(t,r_{k+1}) - \mu r_{k+1}
\sqrt
{n} - \mu Z_{nt}(\omega) \bigr| \geq\frac{\varepsilon}{2} \sqrt
{n} \biggr\}.
\end{eqnarray*}
Taking unions over $r\in[-R,R]$ and $t\in[0,T]$ implies that
\begin{eqnarray*}
&& \Bigl\{ \sup_{r\in[-R,R]} \sup_{t\in[0,T]} \bigl| E_\omega
Y_{n}(t,r) -
\mu r \sqrt{n} - \mu Z_{nt}(\omega) \bigr| \geq\varepsilon\sqrt
{n} \Bigr\} \\
&&\qquad\subset\biggl\{ \max_{k\leq N} \sup_{t\in[0,T]} \bigl| E_\omega
Y_{n}(t,r_k) - \mu r_k \sqrt{n} - \mu Z_{nt}(\omega) \bigr| \geq
\frac{\varepsilon
}{2} \sqrt{n} \biggr\}.
\end{eqnarray*}
Recalling \eqref{finiter} finishes the proof of Theorem \ref{QMCurrent}.
\end{pf*}

%s5 ###
\section{Fluctuations of the centered current} \label{currsec}

Theorems \ref{findimthm} and \ref{findimthmVq} are proved in a similar
way. We spell out
some details for Theorem \ref{findimthm} and restrict to a few remarks
on Theorem \ref{findimthmVq}.
The following representation of the covariance function $\Gamma
((s,q),(t,r))$ will be convenient
(proof by calculus). Recall that $B_\centerdot$ denotes standard
Brownian motion:
%Let $B_\centerdot$ denote a standard Brownian motion with distribution
%$\Pv$.
%Define the function
%$\Gamma$ on $(\R_+\times\R)\times(\R_+\times\R)\times C(\R_+,\R)$ by
%Note that the covariance function $\Gamma((s,q),(t,r))$ has the
%following equivalent definition.
%
%
%e32 ###
\begin{eqnarray}\label{vpGa}
\Gamma((s,q),(t,r))
&=& \mu\int_{-\infty}^{\infty}( \mathbf{P}[B_{\sigma
_1^2s}\le q-x]\mathbf{P}[ B_{\sigma
_1^2t}> r-x] \nonumber\\
&&\hspace*{32pt}{}- \mathbf{P}[B_{\sigma_1^2s}\le q-x, B_{\sigma_1^2t}>
r-x] ) \,dx
\nonumber
\\[-8pt]
\\[-8pt]
\nonumber
&&{} + \sigma_0^2 \biggl\{   \int_{0}^\infty\mathbf
{P}[B_{\sigma_1^2s}\le
q-x]\mathbf{P}[ B_{\sigma_1^2t}\le r-x] \,dx\\
&&\hspace*{32pt}{}   +  \int_{-\infty}^{0} \mathbf{P}[B_{\sigma
_1^2s}> q-x]\mathbf{P}[
B_{\sigma_1^2t}> r-x] \,dx \biggr\}.\nonumber
\end{eqnarray}

Pick time--space points $(t_1,r_1),\ldots,(t_N,r_N)\in\mathbb
R_+\times\mathbb R$ and
$\alpha_1,\ldots, \alpha_N, \beta_1,\ldots,\break \beta_N\in\mathbb
R$. Form the linear combinations
\[
(\bar V_n, \bar Z_n)= \Biggl(n^{-1/4} \sum_{i=1}^N \alpha_i
V_n(t_i,r_i) ,
n^{-1/2} \sum_{i=1}^N \beta_i Z_{nt_i} \Biggr)
\]
and
\[
(\bar V, \bar Z)= \Biggl(  \sum_{i=1}^N \alpha_i V(t_i,r_i) ,
\sum_{i=1}^N \beta_i Z(t_i) \Biggr).
\]
Theorem \ref{findimthm} is proved by showing
$(\bar V_n, \bar Z_n) \stackrel{\mathcal{D}}{\longrightarrow}(\bar V,
\bar Z)$ for an arbitrary choice
of $\{t_i,r_i,\alpha_i,\beta_i\}$.

We can work with $\bar V_n$ alone for a while because much of its analysis
is done under a fixed $\omega$, and then $Z_{n\centerdot}$ is not random:
%
%
%e33 ###
\begin{equation}\qquad\bar V_n=n^{-1/4} \sum_{i=1}^N \alpha_i V_n(t_i,r_i)
=n^{-1/4}\sum_{x\in\mathbb{Z}} \sum_{i=1}^N \alpha_i \bigl[
\mathbf{1}_{ \{ x>0 \} } \phi_{x,i}
 -  \mathbf{1}_{ \{ x\le0 \} } \psi_{x,i}\bigr], \label{vpVnsum}
\end{equation}
where
\[
\phi_{x,i}=\sum_{k=1}^{\eta_0(x)}\mathbf{1}{ \bigl\{ X^{x,k}_{nt_i}\le
nt_i\mathrm{v}_P +r_i\sqrt{n}  \bigr\} }
-E_\omega(\eta_0(x)) P_\omega\bigl\{X^{x}_{nt_i}\le nt_i\mathrm
{v}_P+r_i\sqrt{n} \bigr\}
\]
and
\[
\psi_{x,i}=\sum_{k=1}^{\eta_0(x)}\mathbf{1}{ \bigl\{ X^{x,k}_{nt_i}>
nt_i\mathrm{v}_P+r_i\sqrt{n}  \bigr\} }
-E_\omega(\eta_0(x)) P_\omega\bigl\{X^{x}_{nt_i}> nt_i\mathrm
{v}_P+r_i\sqrt{n} \bigr\}.
\]

Equation \eqref{vpVnsum} expresses $\bar V_n=n^{-1/4}\sum_{x\in
\mathbb Z}
u(x)$ as a sum of
random variables
\[
u(x) = \sum_{i=1}^N \alpha_i \bigl[ \mathbf{1}_{ \{ x>0 \} } \phi_{x,i}
 -  \mathbf{1}_{ \{ x\le0 \} } \psi_{x,i}\bigr]
\]
that are independent and mean zero under the quenched measure $P_\omega
$. They
satisfy
%
%
%e34 ###
\begin{equation}\abs{u(x)}\le\sum_{i=1}^N\,\abs{\alpha_i} \bigl( \eta
_0(x)+E_\omega(\eta
_0(x))\bigr).
\label{vpaux2.6}
\end{equation}

Again we will pick $a(n)\nearrow\infty$ and define
\[
\bar V_n^*=n^{-1/4}\sum_{\abs{x}\le a(n)\sqrt{n}} u(x).
\]
We first show that the rest of the sum can be ignored.

\begin{lem} $ {\lim_{n\to\infty} \mathbb{E}
[(\bar V_n-\bar V_n^*)^2
] =0. }$
\label{vpVV*lm}
\end{lem}
\begin{pf}
By the independence of the $\{u(x)\}$ under $P_\omega$,
%
%
%e35 ###
\begin{eqnarray}\label{vpfd.07}
\hskip12pt&&\mathbb{E}[(\bar V_n-\bar V_n^*)^2]\nonumber\\
&&\qquad = n^{-1/2}E_P\sum
_{\abs{x}> a(n)
\sqrt{n}} E_\omega[u(x)^2]
\nonumber\\
&&\qquad\le Cn^{-1/2}\\
&&\qquad\quad\times\sum_{i=1}^N \sum_{\abs{x}> a(n)\sqrt{n}}
E_P\Biggl[ \mathbf{1}_{ \{ x>0 \} }
\Var_\omega\Biggl(  \sum_{k=1}^{\eta_0(x)}\mathbf{1}{ \bigl\{
X^{x,k}_{nt_i}\le nt_i\mathrm{v}_P +r_i\sqrt{n}  \bigr\} } \Biggr) \nonumber\\
&&\hspace*{90pt}\qquad\quad {}   +   \mathbf{1}_{ \{ x\le0 \} }
\Var_\omega\Biggl(  \sum_{k=1}^{\eta_0(x)}\mathbf{1}{ \bigl\{
X^{x,k}_{nt_i}> nt_i\mathrm{v}_P +r_i\sqrt{n}  \bigr\} }
\Biggr)\Biggr].\nonumber\hskip-12pt
\end{eqnarray}
Consider the first type of variance above:
\begin{eqnarray*}
&&\Var_\omega\Biggl(  \sum_{k=1}^{\eta_0(x)}\mathbf{1}{ \bigl\{
X^{x,k}_{nt_i}\le nt_i\mathrm{v}_P +r_i\sqrt{n}  \bigr\} } \Biggr)\\
&&\qquad = E_\omega(\eta_0(x)) \Var_\omega\bigl(\mathbf{1}{\bigl \{
X^{x}_{nt_i}\le nt_i\mathrm{v}_P +r_i\sqrt{n}  \bigr\} } \bigr)\\
&&\qquad\quad {}+\Var_\omega(\eta_0(x)) P_\omega\bigl\{X^{x}_{nt_i}\le nt_i\mathrm
{v}_P+r_i\sqrt{n} \bigr\}^2 \\
&&\qquad \le[ E_\omega(\eta_0(x)) + \Var_\omega(\eta
_0(x))]
P_\omega\bigl\{X^{x}_{nt_i}\le nt_i\mathrm{v}_P+r_i\sqrt{n} \bigr\}.
\end{eqnarray*}
The upshot is that to show the vanishing of \eqref{vpfd.07} we need to
control terms of the type
%
%
%e36 ###
\begin{equation}\qquad n^{-1/2}\sum_{{x}> a(n)\sqrt{n}} E_P\bigl[
\bigl(E_\omega(\eta(x)) +
\Var_\omega(\eta_0(x))\bigr)
P_\omega\bigl\{X^x_{n}\le n\mathrm{v}_P+r\sqrt{n} \bigr\}\bigr] \label{vpfd.1}
\end{equation}
as $a(n)\to\infty$, together with its counterpart for $x<-a(n)\sqrt{n}$.
For convenience we replaced time
points $nt_i$ with $n$ and $r$ represents $\max r_i$. We treat the part
in~\eqref{vpfd.1} with the variance and omit the rest.
Letting $a_1(n)=a(n)-r$,
\begin{eqnarray*}
&&n^{-1/2}\sum_{{x}> a(n)\sqrt{n}} E_P\bigl[ \Var_\omega(\eta(x))
 P_\omega\bigl\{X^x_{n}\le n\mathrm{v}_P+r\sqrt{n} \bigr\}\bigr] \\
&&\qquad = n^{-1/2}\sum_{{x}> a(n)\sqrt{n}} E_P\bigl[ \Var_\omega
(\eta(x))
 P_{\theta^x\omega}\bigl\{X_{n}\le n\mathrm{v}_P+r\sqrt{n}-x \bigr\}\bigr]
\\
&&\qquad \le n^{-1/2}\sum_{{y}> a_1(n)\sqrt{n}} E_P\bigl[ \Var
_\omega(\eta(0))
 P_{\omega}\bigl\{X_{n}- n\mathrm{v}_P\le-y \bigr\}\bigr] \\
&&\qquad= E_P\biggl[ \Var_\omega(\eta(0))  E_\omega\biggl\{
\biggl( \frac
{X_{n}- n\mathrm{v}_P}{\sqrt{n}}
+{a_1(n)}\biggr)^- \biggr\}  \biggr]\\
&&\qquad\le\{ E_P[(\Var_\omega(\eta(0)))^p]\}^{1/p}
  \biggl\{ E_P\biggl[ \biggl( E_\omega\biggl\{ \biggl( \frac
{X_{n}- n\mathrm{v}_P
}{\sqrt{n}}
+{a_1(n)}\biggr)^- \biggr\} \biggr)^q
\biggr] \biggr\}^{1/q}
\end{eqnarray*}
for some $p>2$ and, hence, $q=p/(p-1)<2$.
By assumption \eqref{vpmomass1}, the first factor above is a constant
if we take $2<p<2+\varepsilon$. Then by the $L^2(\mathbb{P})$
boundedness of
$n^{-1/2} (X_{n}- n\mathrm{v}_P)$ (Proposition~\ref{UIprop}), the second factor
vanishes as
$a(n)\to\infty$.
\end{pf}

Assume now by a truncation that for $\bar V_n^*$ the initial
occupations satisfy
%
%
%e37 ###
\begin{equation}\eta_0(x)\le n^{1/4-\delta} \label{vptrunc}
\end{equation}
for a small $\delta>0$.
Let momentarily $\widetilde V_{n}^*$ denote the variable with truncated
occupations $\tilde\eta_0(x)=\lfloor\eta_0(x)\wedge n^{1/4-\delta}
\rfloor$.

\begin{lem} If
$a(n)\nearrow\infty$ slowly enough,
$\mathbb{E}[\abs{\bar V_{n}^*-\widetilde V_{n}^*}^2 ]\to0$.
\end{lem}
\begin{pf} With\vspace*{-1pt} $A^{x,k}_i$ denoting the random walk events that appear
in $\phi_{x,i}$ and $B^{x,k}_i$ the ones in $\psi_{x,i}$,
\begin{eqnarray*}
&&\bar V_{n}^*-\widetilde V_{n}^* = \sum_{i=1}^N \alpha_i   n^{-1/4}
\Biggl[   \sum_{0<{x}\le a(n)\sqrt{n}} \Biggl(  \sum_{k=\tilde
\eta
_0(x)+1}^{\eta_0(x)}
\mathbf{1}_{ \{ A^{x,k}_i \} }\\
&&\hspace*{141pt}\qquad{} - E_\omega\bigl(\eta_0(x)-\tilde\eta
_0(x)\bigr) P_\omega(A^{x}_i)
\Biggr)\\
&&\hspace*{75pt}\qquad\quad {}- \sum_{ a(n)\sqrt{n}\le x\le0} \Biggl(  \sum_{k=\tilde\eta
_0(x)+1}^{\eta_0(x)}
\mathbf{1}_{ \{ B^{x,k}_i \} }\\
&&\qquad\hspace*{152pt}{} - E_\omega\bigl(\eta_0(x)-\tilde\eta
_0(x)\bigr) P_\omega(B^{x}_i)
\Biggr) \Biggr].
\end{eqnarray*}
Square and
use independence across sites as in the beginning of the proof of Lemma
\ref{vpVV*lm}
to get
\begin{eqnarray*}
E_\omega\abs{\bar V_{n}^*-\widetilde V_{n}^*}^2 &\le& C n^{-1/2}
\sum_{\abs{x}\le a(n)\sqrt{n}} \bigl[ \Var_\omega\bigl( \eta_0(x)-\tilde
\eta
_0(x)\bigr) +
E_\omega\bigl(\eta_0(x)-\tilde\eta_0(x)\bigr)\bigr] \\
&\le& C n^{-1/2} \sum_{\abs{x}\le a(n)\sqrt{n}}
E_\omega\bigl(\eta_0(x)^2\mathbf{1}{ \{ \eta_0(x)\ge n^{1/4-\delta
} \} } \bigr).
\end{eqnarray*}
By shift-invariance,
\[
\mathbb{E}\abs{\bar V_{n}^*-\widetilde V_{n}^*}^2 \le
Ca(n)\mathbb{E}[\eta_0(0)^2\mathbf{1}{ \{ \eta_0(0)\ge
n^{1/4-\delta} \} } ].
\]
Assumption \eqref{vpmomass1} implies that $\mathbb{E}(\eta
_0(0)^2)<\infty$ and, hence,
the last expectation tends to $0$ as $n\to\infty$. The lemma follows.
\end{pf}

Consequently, Theorem \ref{findimthm} is not affected by this truncation.
For the remainder of this proof we work
with the truncated occupation variables that satisfy \eqref{vptrunc}
without indicating it explicitly in
the notation.

Recall that for complex numbers such that $\abs{z_i},\abs{w_i}\le1$,
%
%
%e38 ###
\begin{equation}\Biggl\vert \prod_{i=1}^m z_i  -  \prod_{i=1}^m
w_i
\Biggr\vert
\le\sum_{i=1}^m \abs{z_i-w_i}. \label{auxineq}
\end{equation}

Let
\[
\sigma_{n,\omega}^2(x) = n^{-1/2} E_\omega[u(x)^2].
\]
By \eqref{vpaux2.6} and the truncation \eqref{vptrunc},
%
%
%e39 ###
\begin{equation}\sigma_{n,\omega}^2(x) \le Cn^{-1/2} E_\omega[\eta
_0(x)^2] \le Cn^{-2\delta}
\label{vpaux4}
\end{equation}
which is $<$1 for large enough $n$.
Then
\begin{eqnarray}\label{vpaux5}
&&\biggl\vert E_\omega[e^{i\bar V_{n}^*}] -\prod_{\abs{x}\le a(n)\sqrt{n}}
\biggl(1-\frac12 \sigma_{n,\omega}^2(x)\biggr) \biggr\vert
\nonumber
\\[-8pt]
\\[-8pt]
\nonumber
&&\qquad\le\sum_{\abs{x}\le a(n)\sqrt{n}} \biggl\vert E_\omega\bigl(e^{in^{-1/4}u(x)}\bigr)
-\biggl(1-\frac12 \sigma_{n,\omega}^2(x)\biggr)\biggr\vert
\end{eqnarray}
by an expansion of the exponential, as in the proof of the
Lindeberg--Feller
theorem in \cite{durr}, Section~2.4.b, page~115,
\begin{eqnarray}\label{vpaux7}
&\le&\frac{C\varepsilon(n)}{\sqrt n} \sum_{\abs{x}\le a(n)\sqrt{n}}
E_\omega[u(x)^2]
\nonumber
\\[-8pt]
\\[-8pt]
\nonumber
&&{}+ \frac{C}{\sqrt n} \sum_{\abs{x}\le a(n)\sqrt{n}} E_\omega[u(x)^2
\mathbf{1}{ \{ \abs{u(x)}\ge n^{1/4}\varepsilon(n) \} }   ]
\end{eqnarray}
for some $0<\varepsilon(n)\searrow0$ that we can choose. If
$\varepsilon(n)n^{\delta}\to
\infty$, then
the truncation~\eqref{vptrunc} makes the second sum on line \eqref
{vpaux7} vanish.
Take $E_P$ expectation over the inequalities from \eqref{vpaux5} to
\eqref{vpaux7}.
Since $E_\omega[u(x)^2] \le CE_\omega[\eta_0(x)^2]$,
moment assumption \eqref{vpmomass1} gives
%
%
%e40 ###
\begin{equation}\frac{C\varepsilon(n)}{\sqrt n} \sum_{\abs{x}\le
a(n)\sqrt{n}} \mathbb{E}[u(x)^2]
\le Ca(n)\varepsilon(n). \label{vpaux8}
\end{equation}
Thus, if $a(n)\nearrow\infty$ slowly enough so that $\varepsilon(n)=a(n)
^{-2}\gg n^{-\delta}$,
\eqref{vpaux7} vanishes as $n\to\infty$.

We have reached this intermediate conclusion:
%
%
%e41 ###
\begin{equation}
\lim_{n\to\infty} E_P\biggl\vert E_\omega[e^{i\bar V_{n}^*}]
-\prod_{
\abs{x}\le a(n)\sqrt{n}}
\biggl(1-\frac12 \sigma_{n,\omega}^2(x)\biggr) \biggr\vert=0.
\label{vpaux9}
\end{equation}

The main technical work is encoded in
the following proposition.
Recall the definition
of $\Gamma$ from \eqref{vpGa}.

\begin{prop}\label{vpgprop} There exist bounded continuous functions $g_n$ on
$\mathbb R^N$ with these properties:
\begin{longlist}[(a)]
\item[(a)] $\sup_n \norm{g_n}_\infty<\infty$
and $g_n\to g$ uniformly on compact subsets of $\mathbb R^N$ where $g$
is also
bounded, continuous
and satisfies
%
%
%e42 ###
\begin{eqnarray}\label{vpgGa}
g(z_1,\ldots,z_N)
= \sum_{1\le i,j\le N}\alpha_i\alpha_j \Gamma\bigl((t_i,r_i + z_i),
(t_j,r_j +z_j) \bigr)
\nonumber
\\[-8pt]
\\[-8pt]
\eqntext{\mbox{for $z=(z_1,\ldots,z_N)\in\mathbb R^N$.}}
\end{eqnarray}

\item[(b)] The following limit holds in $P$-probability as $n\to\infty$:
%
%
%e43 ###
\begin{equation}\biggl\vert  \sum_{\abs{x}\le a(n)\sqrt{n}} \sigma
_{n,\omega}^2(x)
-
g_n( n^{-1/2}Z_{nt_1}, \ldots, n^{-1/2}Z_{nt_N} )
\biggr\vert
\longrightarrow0.
\label{vpsigmaglim}
\end{equation}
\end{longlist}
\end{prop}

\begin{pf*}{Proof of Theorem \protect\ref{findimthm} assuming Proposition
\protect\ref{vpgprop}}
By virtue of Lemma \ref{vpVV*lm}, it remains to show
%
%
%e44 ###
\begin{equation}\vert \mathbb{E}[e^{i\bar V_n^*+i\bar Z_n}]
  -  \mathbf{E}[e^{i\bar
V+i\bar Z}] \vert
\to0. \label{vpaux9.2}
\end{equation}
(We need not put coefficients in front of $\bar V_n^*$ and $\bar Z_n$
because these
coefficients can be subsumed in the $\alpha_i,\beta_i$ coefficients.)
Define the random $N$-vectors
\[
\mathbf{z}_n^{1,N}=(n^{-1/2}Z_{nt_1},\ldots, n^{-1/2}Z_{nt_N})
\quad  \mbox{and}\quad  \mathbf{z}^{1,N}=(Z(t_1),\ldots, Z(t_N)).
\]
Then the conditional
distribution of $V$ given $Z$, described in conjunction with
\eqref{vpcov} above, together with \eqref{vpgGa} gives
\[
\mathbf{E}[e^{i\bar V+i\bar Z}] = \mathbf{E}\bigl[e^{ -{(1/2)} g(\mathbf
{z}^{1,N}) +i\bar Z}\bigr].
\]
Now bound the absolute value in \eqref{vpaux9.2} by
\begin{eqnarray}\label{vpaux9.25}
&&\bigl\vert E_P\bigl[ E_\omega(e^{i\bar V_n^*})   e^{i\bar
Z_n(\omega)}
\bigr]   -
\mathbf{E}\bigl[e^{-(1/2) g(\mathbf{z}^{1,N}) +i\bar Z}\bigr] \bigr\vert
\nonumber\\
&&\qquad\le E_P\bigl\vert E_\omega(e^{i\bar V_{n}^*}) - e^{ -(1/2)
g_n(\mathbf{z}
_n^{1,N}) }  \bigr\vert\\
 &&\quad\qquad{} +  \bigl \vert E_P\bigl[e^{-(1/2) g_n(\mathbf{z}_n^{1,N}) +i\bar
Z_n(\omega
)}\bigr]   -
\mathbf{E}\bigl[e^{-(1/2) g(\mathbf{z}^{1,N}) +i\bar Z}\bigr]
\bigr\vert.\nonumber
\end{eqnarray}
The last absolute values expression above vanishes as $n\to\infty$ by the
invariance principle $n^{-1/2}Z_{n\cdot}\stackrel{\mathcal
{D}}{\longrightarrow}Z(\cdot)$ [Theorem \ref
{QCLTthm}, part (2)]
and by a simple property of weak convergence stated in Lemma \ref{vpweakcomvlm}
after this proof.
The second-to-last term
is bounded as follows:
\begin{eqnarray}
&& E_P\bigl\vert E_\omega(e^{i\bar V_{n}^*}) - e^{ -(1/2)
g_n(\mathbf{z}
_n^{1,N}) }  \bigr\vert\nonumber\\
  &&\qquad\le
E_P\biggl\vert E_\omega(e^{i\bar V_{n}^*}) -\prod_{\abs{x}\le a(n)\sqrt{n}}
\biggl(1-\frac12 \sigma_{n,\omega}^2(x)\biggr) \biggr\vert\label
{vpaux9.3} \\
&&\qquad\quad {} + E_P\biggl\vert  \prod_{\abs{x}\le a(n)\sqrt{n}}
\biggl(1-\frac12 \sigma_{n,\omega}^2(x)\biggr)  -
\exp\biggl\{ -\frac12 \sum_{\abs{x}\le a(n)\sqrt{n}} \sigma_{n,\omega}^2(x)
\biggr\}
 \biggr\vert\label{vpaux9.4}\hskip-12pt \\
&&\qquad\quad {} + E_P\biggl\vert
  \exp\biggl\{ -\frac12 \sum_{\abs{x}\le a(n)\sqrt{n}} \sigma
_{n,\omega}^2(x)
\biggr\}
- \exp\biggl( -\frac12 g_n(\mathbf{z}_n^{1,N}) \biggr)
\biggr\vert. \label{vpaux9.5}
\end{eqnarray}
Let $n\to\infty$.
Line \eqref{vpaux9.3} after the inequality vanishes by \eqref{vpaux9}.
Line \eqref{vpaux9.4} vanishes by the inequalities
\begin{eqnarray*}
\exp\biggl(  -\frac12(1+n^{-2\delta})      \sum_{\abs{x}\le a(n)
\sqrt{n}} \sigma_{n,\omega}^2(x) \biggr)
&\le&\prod_{\abs{x}\le a(n)\sqrt{n}}
     \biggl(1-\frac12 \sigma_{n,\omega}^2(x)\biggr)\\
&\le&\exp\biggl(  -\frac12\sum_{\abs{x}\le a(n)\sqrt{n}} \sigma
_{n,\omega
}^2(x) \biggr),
\end{eqnarray*}
where we used \eqref{vpaux4} and $ -y-y^2\le\log(1-y)\le-y$ for
small $y>0$.
Finally, line \eqref{vpaux9.5} vanishes by \eqref{vpsigmaglim}.

We have shown that line \eqref{vpaux9.25} vanishes as $n\to\infty$
and thereby
verified \eqref{vpaux9.2}.
This completes the proof of Theorem \ref{findimthm}, assuming
Proposition \ref{vpgprop}.
\end{pf*}

Lines \eqref{vpaux9.3}--\eqref{vpaux9.5}, $\mathbf{z}_n^{1,N}\mathop{\stackrel{\mathcal{D}}{\longrightarrow}}\limits_{n\to
\infty}\mathbf{z}^{1,N}$
and $g_n\to g$ uniformly on compacts show that
%
%
%e45 ###
\begin{equation}
E_P\bigl\vert E_\omega(e^{i\bar V_{n}^*}) - e^{ -(1/2) g(\mathbf
{z}^{1,N}) }
 \bigr\vert\to0.
\label{vplim7}
\end{equation}
This verifies the remark stated after Theorem \ref{findimthm}.%\vadjust{\goodbreak}

The next lemma was used in the proof above. We omit its short and
simple proof.\vspace*{-3pt}

\begin{lem}
Suppose $\zeta_n\stackrel{\mathcal{D}}{\longrightarrow}\zeta$ for
random variables with values in some
Polish space $S$.
Let
$f_n,f$ be bounded, continuous functions on $S$ such that $\sup_n
\norm{f_n}_\infty<\infty$
and $f_n\to f$
uniformly on compact sets. Then $f_n(\zeta_n)\stackrel{\mathcal
{D}}{\longrightarrow}f(\zeta)$.
\label{vpweakcomvlm}\vspace*{-3pt}
\end{lem}

We turn to the proof of the main
technical proposition, Proposition \ref{vpgprop}.\vspace*{-3pt}

\begin{pf*}{Proof of Proposition \protect\ref{vpgprop}}
Consider $n$ large enough so that $a(n)>\max_i\abs{r_i}$:
%
%
%e46 ###
\begin{eqnarray}\label{vpfindim1}
\sum_{\abs{x}\le a(n)\sqrt{n}} \sigma_{n,\omega}^2(x) &=& n^{-1/2}\sum_{
\abs{x}\le a(n)\sqrt{n}} E_\omega[u(x)^2]\nonumber\\
&=& n^{-1/2}\sum_{\abs{x}\le a(n)\sqrt{n}} \Cov_\omega[u(x),u(x)]
\nonumber
\\[-8.5pt]
\\[-8.5pt]
\nonumber
&=&\sum_{1\le i,j\le N}\alpha_i\alpha_j   n^{-1/2}\sum_{\abs{x}\le a(n)
\sqrt{n}}
\bigl[ \mathbf{1}_{ \{ x>0 \} } \Cov_\omega(\phi_{x,i}, \phi_{x,j})\\
  &&\hspace*{130pt}{}+  \mathbf{1}_{ \{ x\le0 \} } \Cov_\omega(\psi_{x,i}, \psi
_{x,j})\bigr].\nonumber
\end{eqnarray}
Whenever we work with a fixed $(i,j)$ we let $((s,q),(t,r))$ represent
$((t_i,r_i),\break(t_j, r_j))$
to avoid
excessive subscripts.
To each term above apply the formula for the covariance of two random
sums, with
$\{Z_i\}$ i.i.d. and independent of~$K$:
\[
\Cov\Biggl( \sum_{i=1}^K f(Z_i) , \sum_{j=1}^K g(Z_j)\Biggr) =
EK  \Cov(f(Z),g(Z))+ \Var(K) Ef(Z)  Eg(Z).
\]
The first covariance on the last line of \eqref{vpfindim1} develops as
%
%
%e47 ###
\begin{eqnarray}%
\label{vpfindim3}
&&\Cov_\omega(\phi_{x,i}, \phi_{x,j})\nonumber\\
&&\qquad= E_\omega(\eta_0(x))
P_\omega\bigl\{X^{x}_{ns}\le ns\mathrm{v}_P+q\sqrt{n},  X^{x}_{nt}\le
nt\mathrm{v}_P+r\sqrt{n}
\bigr\} \nonumber\\
&&\qquad\quad {} - E_\omega(\eta_0(x))
P_\omega\bigl\{X^{x}_{ns}\le ns\mathrm{v}_P+q\sqrt{n} \bigr\}P_\omega\bigl\{
X^{x}_{nt}\le nt\mathrm{v}_P+r\sqrt
{n}  \bigr\} \nonumber\\
&&\qquad\quad {} + \Var_\omega(\eta_0(x))
P_\omega\bigl\{X^{x}_{ns}\le ns\mathrm{v}_P+q\sqrt{n} \bigr\}P_\omega\bigl\{
X^{x}_{nt}\le nt\mathrm{v}_P+r\sqrt
{n}  \bigr\} \\
&&\qquad= - E_\omega(\eta_0(x))
P_\omega\bigl\{X^{x}_{ns}\le ns\mathrm{v}_P+q\sqrt{n},  X^{x}_{nt}>
nt\mathrm{v}_P+r\sqrt{n}  \bigr\}
\nonumber\\
&&\qquad\quad {} +  E_\omega(\eta_0(x))
P_\omega\bigl\{X^{x}_{ns}\le ns\mathrm{v}_P+q\sqrt{n} \bigr\}P_\omega\bigl\{
X^{x}_{nt}> nt\mathrm{v}_P+r\sqrt
{n}  \bigr\} \nonumber\\
&&\qquad\quad {}  + \Var_\omega(\eta_0(x))
P_\omega\bigl\{X^{x}_{ns}\le ns\mathrm{v}_P+q\sqrt{n} \bigr\}P_\omega\bigl\{
X^{x}_{nt}\le nt\mathrm{v}_P+r\sqrt
{n}  \bigr\}.\nonumber
\end{eqnarray}

Develop the second covariance in a similar vein, and then collect the terms:
\begin{eqnarray}\label{vpfindim3.3}
&&\hspace*{-4pt}\sum_{\abs{x}\le a(n)\sqrt{n}} \sigma_{n,\omega}^2(x) \nonumber\\
&&\hspace*{-4pt}\qquad= \sum_{1\le i,j\le N}\alpha_i\alpha_j  \biggl[
n^{-1/2}\sum_{\abs{x}\le a(n)\sqrt{n}} E_\omega(\eta_0(x))\bigl(
P_\omega\bigl\{X^{x}_{nt_i}\le nt_i\mathrm{v}_P+r_i\sqrt{n} \bigr\}
\nonumber\\
\label{vpfindim3.4}%&&\hspace*{-3pt}\qquad\hspace*{132pt}\quad{}\times\bigl(
%P_\omega\bigl\{X^{x}_{nt_i}\le nt_i\mathrm{v}_P+r_i\sqrt{n} \bigr\}\nonumber\\
&&\qquad\quad\hspace*{179pt} {}\times P_\omega
\bigl\{ X^{x}_{nt_j}>
nt_j\mathrm{v}_P+r_j\sqrt{n}  \bigr\}\\
&&\qquad\quad \hspace*{179pt}{}-  P_\omega\bigl\{X^{x}_{nt_i}\le
nt_i\mathrm{v}_P+r_i\sqrt{n},\nonumber
\\[-8pt]
\\[-8pt]
\label{vpfindim3.5}&&\qquad\quad \hspace*{209pt} X^{x}_{nt_j}> nt_j\mathrm{v}_P
+r_j\sqrt{n}  \bigr\}\bigr)
\nonumber
\\[-8pt]
\\[-8pt]
\label{vpfindim3.8}&&\qquad\hspace*{59pt}\quad {}+n^{-1/2}\sum_{\abs{x}\le a(n)\sqrt{n}} \Var_\omega(\eta_0(x)) \nonumber\\
&&\qquad\hspace*{151pt}{}\times \bigl(
\mathbf{1}_{ \{ x>0 \} }
P_\omega\bigl\{X^{x}_{nt_i}\le nt_i\mathrm{v}_P+r_i\sqrt{n} \bigr\}\nonumber\\
&&\qquad\quad \hspace*{154pt}{} \times P_\omega
\{ X^{x}_{nt_j}\le
nt_j\mathrm{v}_P+r_j\sqrt{n}  \} \\
&&\qquad\quad \hspace*{154pt}{}   + \mathbf{1}_{ \{ x\le0 \} }
P_\omega\bigl\{X^{x}_{nt_i}> nt_i\mathrm{v}_P+r_i\sqrt{n} \bigr\}\nonumber\\
&&\qquad\quad \hspace*{177pt}{} \times P_\omega\bigl\{
X^{x}_{nt_j}> nt_j\mathrm{v}_P
+r_j\sqrt{n}  \bigr\}  \bigr) \biggr]. \nonumber
\end{eqnarray}

The function $g_n(z_1,\ldots,z_N)$ required for Proposition \ref{vpgprop}
is defined as the linear combination of integrals
of Brownian probabilities\vspace*{1pt} that match up with the terms of the sum above.
For $(z_1,\ldots,z_N)\in\mathbb R^N$,
%
%
%e48 ###
\begin{eqnarray}\label{vpg_n}
&&\hspace*{-4pt}g_n(z_1,\ldots,z_N)\nonumber\\
&&\qquad= \sum_{1\le i,j\le N}\alpha_i\alpha_j
\biggl[   \mu\int_{-a(n)}^{a(n)}( \mathbf{P}[B_{\sigma
_1^2t_i}\le
z_i+r_i-x]\nonumber\\
&&\qquad\hspace*{113pt} {}\times \mathbf{P}[ B_{\sigma_1^2t_j}> z_j+r_j-x]\nonumber\\
&&\qquad{}\hspace*{113pt}
-  \mathbf{P}[B_{\sigma_1^2t_i}\le z_i+r_i-x, B_{\sigma_1^2t_j}\nonumber\\
&&\qquad{}\hspace*{183pt}>
z_j+r_j-x] )
\,dx \\
& &\hspace*{93pt}{}  + \sigma_0^2 \biggl\{   \int_{0}^{a(n)}
\mathbf{P}[B_{\sigma
_1^2t_i}\le z_i+r_i-x]\nonumber\\
&&\hspace*{148pt}{}\times\mathbf{P}[ B_{\sigma_1^2t_j}\le z_j+r_j-x]
\,dx\nonumber\\
&&\hspace*{122pt}{}   +  \int_{-a(n)}^{0}
\mathbf{P}[B_{\sigma_1^2t_i}> z_i+r_i-x]\nonumber\\
&&\hspace*{162pt}{}\times \mathbf{P}[ B_{\sigma
_1^2t_j}> z_j+r_j-x] \,dx
\biggr\}   \biggr].\nonumber
\end{eqnarray}
Let $g(z_1,\ldots,z_N)$ be the function defined by the above sum of
integrals with
$a(n)$ replaced by $\infty$. Then \eqref{vpgGa} holds by direct
comparison with
definition \eqref{vpGa}. Part~(a) of Proposition \ref{vpgprop} is now clear.

%The expression in brackets on lines\eqref{vpfindim3.2}--
%probabilities $\{\w_x\}$.

To prove limit \eqref{vpsigmaglim} in part (b) of Proposition \ref
{vpgprop}, namely, that
\[
\biggl\vert  \sum_{\abs{x}\le a(n)\sqrt{n}} \sigma_{n,\omega}^2(x)
  -
g_n( n^{-1/2}Z_{nt_1}, \ldots, n^{-1/2}Z_{nt_N} )
\biggr\vert
\stackrel{P}\longrightarrow0,
\]
we approximate the
sums on lines \eqref{vpfindim3.3}--\eqref{vpfindim3.8} with the corresponding
integrals from \eqref{vpg_n}.
The steps are the same for each sum. We illustrate this reasoning with
the sum of the terms on
line \eqref{vpfindim3.4}, given by
%
%
%e49 ###
\begin{eqnarray}\label{vpdefU}
&& U_n(\omega)=\sum_{\abs{m}\le a(n)\sqrt{n}} E_\omega(\eta_0(m))
\nonumber
\\[-8pt]
\\[-8pt]
\nonumber
&&\hspace*{64pt}\qquad{}\times P_\omega\bigl\{X^{m}_{ns}\le ns\mathrm{v}_P+q\sqrt{n},  X^{m}_{nt}>
nt\mathrm{v}_P+r\sqrt{n}  \bigr\}
\end{eqnarray}
and the corresponding part of \eqref{vpg_n}, defined by
%
%
%e50 ###
\begin{eqnarray}\label{vpdefU*}
&& U^*_n(\omega)=\mu\int_{-a(n)}^{a(n)} \mathbf{P}\biggl\{B_{\sigma
_1^2s}\le\frac
{Z_{ns}(\omega)}{\sqrt{n}}-x+q ,
\nonumber
\\[-8pt]
\\[-8pt]
\nonumber
&&\qquad\hspace*{77pt}{} B_{\sigma_1^2t}> \frac{Z_{nt}(\omega)}{\sqrt{n}}-x+r \biggr\}\, dx.
\end{eqnarray}

The goal is to show
\[
\lim_{n\to\infty}\abs{ n^{-1/2} U_n(\omega)- U^*_n(\omega)}=0\qquad
\mbox{in $P$-probability. }
\]
The steps are the same as those employed in the proofs of
Lemmas \ref{WtildeW}--\ref{hWbW}. First
approximate $U_n$ with
%
%
%e51 ###
\begin{eqnarray}\label{vpdefwtU}
&&\widetilde U_n(\omega)=\sum_{\abs{m}\le a(n)\sqrt{n}} E_\omega(\eta_0(m))
\mathbf{P}\biggl\{B_{\sigma_1^2s}\le\frac{Z_{ns}(\theta^m\omega
)}{\sqrt{n}}-\frac
{m}{\sqrt{n}}+q ,
\nonumber
\\[-8pt]
\\[-8pt]
\nonumber
&&\hspace*{133pt}\qquad{} B_{\sigma_1^2t}> \frac{Z_{nt}(\theta^m\omega)}{\sqrt{n}}-\frac
{m}{\sqrt{n}}+r
\biggr\}.
\end{eqnarray}
This approximation is similar to the proof of Lemma \ref{WtildeW} and
uses the fact that, for a fixed $s,t>0$, the limits of the form \eqref
{QCLT} are uniform in $x,y\in\mathbb R$.
Then remove the shift from $Z_n(\omega)$ by defining
%
%
%e52 ###
\begin{eqnarray}\label{vpdefwhU}
&&\widehat U_n(\omega)=\sum_{\abs{m}\le a(n)\sqrt{n}} E_\omega(\eta_0(m))
\mathbf{P}\biggl\{B_{\sigma_1^2s}\le\frac{Z_{ns}(\omega)}{\sqrt
{n}}-\frac{m}{\sqrt
{n}}+q ,
\nonumber
\\[-8pt]
\\[-8pt]
\nonumber
&&\qquad\hspace*{132pt}{} B_{\sigma_1^2t}> \frac{Z_{nt}(\omega)}{\sqrt{n}}-\frac{m}{\sqrt
{n}}+r \biggr\}
\end{eqnarray}
and showing that
$ \lim_{n\to\infty}n^{-1/2} \abs{\widetilde U_n-\widehat
U_n}=0$, in $P$-probability.
For the last step, to show
$ \lim_{n\to\infty}\abs{ n^{-1/2} \widehat U_n(\omega)-
U^*_n(\omega)}=0$ in
$P$-probability,
truncate the sum~\eqref{vpdefwhU} and the integral \eqref{vpdefU*},
use a Riemann approximation of the sum, introduce an intermediate scale
for further partitioning and appeal to the ergodic theorem, as was done
in Lemma \ref{hWbW}. We omit these details
since the corresponding steps were spelled out in full in Section \ref
{qmeansec}.

We have verified the part of the desired limit \eqref{vpsigmaglim} that
comes from pairing up
the sum on line \eqref{vpfindim3.4} with the second line of \eqref{vpg_n}.
The remaining parts are handled similarly.
This completes the proof of Proposition \ref{vpgprop}.
\end{pf*}

Theorem \ref{findimthm} has now been proved.
Proof of Theorem \ref{findimthmVq} goes essentially the same way.
% as the proof of Theorem \ref{findimthm}.
The crucial difference comes at the point \eqref{vpdefwhU} where
$\widehat U_n$ is introduced. Instead of $n^{-1/2}Z_{ns}(\omega)$ and
$n^{-1/2}Z_{nt}(\omega)$
inside the Brownian probability $\mathbf{P}$, one has
$n^{-1/2}(Z_{ns}(\omega)-Z_{ns}(\theta^m\omega))$ and
$n^{-1/2}(Z_{nt}(\omega)-Z_{nt}(\theta^m\omega))$.
These vanish on the scale considered here, with $\abs{m}\le a(n)\sqrt n$,
by the arguments used in the proof of Lemma \ref{tWhW}.

Consequently, in the subsequent approximation by $U_n^*$ at \eqref{vpdefU*},
the terms $n^{-1/2}Z_{ns}(\omega)$ and $n^{-1/2}Z_{nt}(\omega)$ have
disappeared.
Then in limit \eqref{vpaux9.2} in Proposition \ref{vpgprop} we can take
$g_n(0,\ldots,0)$.

\begin{appendix}
%s6 ###
\section*{Appendix: Uniform integrability of $\sup_{k\leq n} (X_k-k\mathrm{v}_P)/\sqrt{n}$}\label{UIapp}

In this appendix we give the proof of Proposition \ref{UIprop}. The
main tool used in the proof is a martingale representation that was
given in the proof of the averaged central limit theorem in \cite{zRWRE}.
% and a moment bound given in \cite{gQCLT} that helps control the
%fluctuations of the function $h(x,\w)$ defined in \eqref{hdef}.
%
%where $\s_1^2 = \vp^3 E_P \Var_\w T_1$ and $\s_2^2 = \vp^2 \Var( E_\w
%T_1 )$.
%Moreover, there exists a constant $C<\infty$ such that
% \E[\sup_{k\leq n} (X_k - k \vp)^2 ] \leq C n.
%
%We will use the martingale representation used in proving the averaged
%central limit theorem in \cite{zRWRE}. For any environment $\w$, let
% h(x,\w) :=
%0 & x= 0 \\
%- \vp\sum_{i=x}^{-1} ( E_{\theta^i \w} T_1 - \E T_1 ) & x \leq-1.
%Also, let $\mathcal{F}_n := \s(X_i : i\leq n)$.
Recall the definition of $h(x,\omega)$ in \eqref{hdef}, and let
$\mathcal
{F}_n := \sigma(X_i \dvtx  i\leq n)$.
Then, $M_n := X_n - n\mathrm{v}_P+ h(X_n, \omega)$ is an $\mathcal
{F}_n$-martigale
under the measure $P_\omega$.
The correction term $h(X_n,\omega)$ may further be decomposed as
$h(X_n,\omega)
= Z_n(\omega) + R_n$, where $Z_n(\omega) = h(\lfloor n\mathrm{v}_P
\rfloor,\omega)$ and $R_n := h(X_n,\omega
)- Z_n(\omega) $.
The main contributions to $X_n - n\mathrm{v}_P$ come from $M_n$ and
$Z_n(\omega)$,
while the term $R_n$ contributes on a scale of order less than $\sqrt
{n}$. $M_n$ accounts for the fluctuations due to the randomness of the
walk in a fixed environment, and $Z_n(\omega)$ accounts for the
fluctuations due to randomness of the environment.

Using the above notation, we then have
%
%
%e53 ###
\setcounter{equation}{0}
\begin{eqnarray}\label{Mnexpansion}
\mathbb{E}( X_{n} - n\mathrm{v}_P)^2 &=& \mathbb{E}M_n^2
+ E_P Z_n(\omega)^2 + \mathbb{E}R_n^2
\nonumber
\\[-8pt]
\\[-8pt]
\nonumber
&&{} - 2 \mathbb{E}[ M_n R_n ]+ 2 \mathbb{E}[ Z_n(\omega) R_n].
\end{eqnarray}
Note that the term $\mathbb{E}[ M_n Z_n(\omega) ]$ is missing on the
right-hand
side above. This is because $Z_n(\omega)$ depends only on the environment
and $M_n$ is a martingale under~$P_\omega$ and, thus, $\mathbb{E}[M_n
Z_n(\omega)] =
E_P[ Z_n(\omega) E_\omega(M_n) ] = 0$.
%Since $Z_n(\w)$ depends only on the environment and $M_n$ is a
%martingale under $P_\w$, $\E[M_n Z_n(\w)] = E_P[ Z_n(\w) E_\w
%(M_n) ] = 0$.
%This eliminates the fourth term in \eqref{Mnexpansion}.
Since H\"{o}lder's inequality implies that
\[
\mathbb{E}[ M_n R_n ] + \mathbb{E}[ Z_n(\omega) R_n] \leq\bigl(
(\mathbb{E}M_n^2)^{1/2} + (E_P
Z_n(\omega)^2)^{1/2} \bigr)(\mathbb{E}R_n^2)^{1/2},
\]
to complete the proof of \eqref{UIclaim}, it is enough to show
%
%
%e54 ###
\begin{equation} \label{MnZnL2}
\lim_{n\rightarrow\infty} \frac{1}{n} \mathbb{E}M_n^2 = \sigma_1^2,\qquad
\lim_{n\rightarrow\infty} \frac{1}{n} E_P Z_n(\omega)^2 = \sigma_2^2
\end{equation}
and
%
%
%e55 ###
\begin{equation} \label{RnL2}
\lim_{n\rightarrow\infty} \frac{1}{n} \mathbb{E}R_n^2 = 0.
\end{equation}
Since $Z_n(\omega) = h(\lfloor n\mathrm{v}_P \rfloor,\omega)$, to
prove the second statement in
\eqref{MnZnL2}, it is enough to show that
\[
\lim_{n\rightarrow\infty} \frac{1}{n} E_P[ h(n,\omega)^2 ] =
\mathrm{v}_P\Var( E_\omega T_1) =
\frac{1}{\mathrm{v}_P}\sigma_2^2.
\]
However, since $h(n,\omega)$ is the sum of mean zero terms,
\begin{eqnarray*}
E_P[h(n,\omega)^2] &=& \Var(h(n,\omega)) \nonumber\\
&=& \mathrm{v}_P^2 \sum
_{i=0}^{n-1} \Var( E_\omega T_1
) + 2 \mathrm{v}_P^2 \sum_{0\leq i < j \leq n-1} \Cov(
E_{\theta^i \omega} T_1,
E_{\theta^j \omega} T_1 ) \\
&=& n \mathrm{v}_P^2 \Var( E_\omega T_1 ) + 2 \mathrm{v}_P^2 \sum
_{k=1}^{n-1} (n-k) \Cov(
E_{\omega} T_1, E_{\theta^k \omega} T_1 ),\nonumber
\end{eqnarray*}
where the last equality is due to the shift invariance of environments.
Since $E_{\theta^k \omega} T_1 = 1 + \rho_k + \rho_k E_{\theta
^{k-1} \omega}
T_1$ (see the derivation of a formula for $E_\omega T_1$ in \cite
{sRWRE} or
\cite{zRWRE}), the fact that $P$ is an i.i.d. law on environments
implies that
\[
\Cov( E_{\omega} T_1, E_{\theta^k \omega} T_1 ) =
(E_P\rho_0) \Cov
( E_{\omega} T_1, E_{\theta^{k-1} \omega} T_1 ).
\]
Iterating this computation, we get that $\Cov( E_{\omega} T_1,
E_{\theta^k \omega} T_1 ) = (E_P \rho_0)^k \Var( E_\omega
T_1 )$.
Therefore,
\begin{eqnarray*}
E[h(n,\omega)^2] &=& n \mathrm{v}_P^2 \Var( E_\omega T_1 ) + 2
\mathrm{v}_P^2 \Var( E_\omega T_1) \sum
_{k=1}^{n-1} (n-k) (E_P \rho_0)^k \\
& =& n \mathrm{v}_P^2 \Var( E_\omega T_1 ) \Biggl( 1 + 2 \sum
_{k=1}^{n-1} (E_P \rho
_0)^k \Biggr)\\
&&{} - 2\mathrm{v}_P^2 \Var(E_\omega T_1) \sum_{k=1}^{n-1}
k (E_P \rho_0)^k .
\end{eqnarray*}
Since $E_P \rho_0 < 1$, this implies that
%
%
%e56 ###
\begin{eqnarray}\label{hL2}
\lim_{n\rightarrow\infty} \frac{1}{n} E_P[h(n,\omega)^2] &= &\mathrm
{v}_P^2 \Var( E_\omega T_1
)\biggl( 1 + 2 \frac{ E_P \rho_0 }{1- E_P \rho_0} \biggr)
\nonumber
\\[-8pt]
\\[-8pt]
\nonumber
&=& \mathrm
{v}_P\Var(
E_\omega T_1 ),
\end{eqnarray}
where the last equality is from the explicit formula for $\mathrm
{v}_P$ given in
\eqref{LLN}.
Thus, we have proved the second statement in \eqref{MnZnL2}.

We now turn to the proof of the first statement in \eqref{MnZnL2}. Let
\[
V_n := \sum_{k=1}^n E_\omega[ (M_{k+1} - M_k)^2 | \mathcal
{F}_k
] .
\]
Note that $E_\omega V_n = E_\omega M_n^2$ since $M_n$ is a martingale under
$P_\omega$. Thus, the first statement in \eqref{MnZnL2} is equivalent to
$\lim_{n\rightarrow\infty} \mathbb{E}V_n /n = \sigma_1^2$.
A direct computation (see the proof of the averaged central limit
theorem on page 211 of \cite{zRWRE}) yields that $E_\omega[
(M_{k+1} -
M_k)^2 | \mathcal{F}_k ] = g(\theta^{X_k} \omega)$,
where
\[
g(\omega) = \mathrm{v}_P^2 \bigl( \omega_0 (E_\omega T_1 - 1)^2 +
(1-\omega_0)(E_{\theta^{-1}\omega}
T_1 + 1)^2 \bigr).
\]
%
%which is the environment viewed from the particles current location.
Recall the definition of $f(\omega)$ in \eqref{fdef}, and let $Q$ be a
measure on environments defined by $\frac{dQ}{dP}(\omega) = f(\omega
)$, where
$f(\omega)$ is defined in \eqref{fdef}.
Under the averaged measure $\mathbb Q( \cdot) = E_Q[ P_\omega(\cdot)
] $, the
sequence $\{ \theta^{X_k} \omega\}_{k\in\mathbb N}$ is stationary
and ergodic.
Therefore, $\frac{V_n}{n} = \frac{1}{n} \sum_{k=1}^n g(\theta^{X_k}
\omega
)$ converges in $L^1(\mathbb Q)$ to
\[
E_Q[ g(\omega)] = E_P\biggl[ \frac{dQ}{dP}(\omega) g(\omega)
\biggr] = \mathrm{v}_P^3 E_P[
\Var_\omega T_1 ] = \sigma_1^2,
\]
where the second to last equality follows from the formulas for $\frac
{dQ}{dP}(\omega)$ and $g(\omega)$ given above, the explicit formula
for $\Var
_\omega T_1$ shown in \cite{pThesis}, and the shift invariance of the
law $P$.
Since $\frac{dQ}{dP}(\omega) = f(\omega) \geq\mathrm{v}_P$, we
obtain that
\[
\mathbb{E}| V_n/n - \sigma_1^2 | = E_Q \biggl[ \frac
{dP}{dQ}(\omega) E_\omega
| V_n/n - \sigma_1^2 | \biggr] \leq\frac{1}{\mathrm
{v}_P} E_{\mathbb Q}|
V_n/n- \sigma_1^2 | \mathop{\longrightarrow}_{n\rightarrow\infty
}0.
\]
Thus, since $V_n/n$ converges in $L^1(\mathbb Q)$ to $\sigma_1^2$,
$V_n/n$ also
converges to $\sigma_1^2$ in $L^1(\mathbb{P})$.

Finally, we turn to the proof of \eqref{RnL2}.
Fix a $\beta\in(1/2,1)$. Since $R_n = h(X_n, \omega) - h(\lfloor
n\mathrm{v}_P \rfloor, \omega)$,
%we have for any $\b\in(1/2,1)$ that
%
\begin{eqnarray*}
E_\omega R_n^2 &\leq&\sup_{x\dvtx  |x-\lfloor n\mathrm{v}_P \rfloor| \leq
n^{\beta}} |h(x,\omega) - h(\lfloor n\mathrm{v}_P \rfloor,\omega
)|^2 \\
&&{}+ \sup_{|x|\leq n} 4|h(x,\omega)|^2 P_\omega(|X_n - \lfloor
n\mathrm{v}_P \rfloor| >
n^{\beta}).
\end{eqnarray*}
Then, the shift invariance of the measure $P$ and H\"older's inequality
imply that, for any $\delta>0$,
\begin{eqnarray*}
\mathbb{E}R_n^2 & \leq&2E_P \Bigl[ \sup_{|x|\leq n^{\beta}}
h(x,\omega)^2 \Bigr] + 4
E_P\Bigl[ \sup_{|x|\leq n} h(x,\omega)^2 P_\omega(|X_n - \lfloor
n\mathrm{v}_P \rfloor| > n^{\beta})
\Bigr] \\
%&\leq( E_P [ \sup_{|x|\leq n^{\b}} h(x,\w)^{2+2\d} ]
%)^{1/(1+\d)} + 4 (E_P[ \sup_{|x|\leq n} h(x,\w)^{2+2
&\leq& E_P \Bigl[ \sup_{|x|\leq n^{\beta}} h(x,\omega)^2 \Bigr] +4\Bigl(E_P\Bigl[ \sup_{|x|\leq n} h(x,\omega)^{2+2\delta} \Bigr]
\Bigr)^{1/(1+\delta)} \\
&&\hspace*{98pt}{}\times\mathbb{P}
(|X_n - \lfloor n\mathrm{v}_P \rfloor| > n^{\beta})^{\delta
/(1+\delta)} \\
&\leq& C n^{\beta} + C n \mathbb{P}(|X_n - \lfloor n\mathrm{v}_P
\rfloor| > n^{\beta})^{\delta/(1+\delta)},
\end{eqnarray*}
where the last inequality follows from Lemma \ref{Goldsheid}.
%Let $h^*(n,\w) := \sup_{|x|\leq n} |h(x,\w)|$. Then, the shift
%invariance of the measure $P$ and H\"older's inequality imply that for
%any $\d>0$,
% \E R_n^2 & \leq E_P [ h^*(n^{\b},\w)^2 ] + 4 E_P[
%h^*(n,\w)^2 P_\w(|X_n - \fl{n\vp}| > n^{\b}) ] \\
%&\leq( E_P [ h^*(n^{\b},\w)^{2+2\d} ] )^{1/(1+
%Goldsheid has shown \cite[Lemma 4]{gQCLT} that for any $\d>0$
%sufficiently small, there exists a constant $C<\infty$ such that $
%(E_P[ h^*(n,\w)^{2+2\d}])^{1/(1+\d)} \leq C n$.
%Thus,
The first term on the right above is $o(n)$ since $\beta< 1$, and the
second term on the right is $o(n)$ because $\beta>1/2$ and the averaged
central limit theorem implies that $\mathbb{P}(|X_n - \lfloor n\mathrm
{v}_P \rfloor| > n^{\beta})$
tends to zero. This completes the proof of \eqref{RnL2} and thus also
the first part of Proposition~\ref{UIprop}.

To prove the second part of Proposition \ref{UIprop}, we again use the
representation $X_n - n\mathrm{v}_P= M_n - Z_n(\omega) - R_n$. Then,
\[
\mathbb{E}\Bigl[ \sup_{k\leq n} (X_k - k\mathrm{v}_P)^2 \Bigr]
\leq3\mathbb{E}\Bigl[ \sup
_{k\leq n} M_k^2 \Bigr] + 3\mathbb{E}\Bigl[ \sup_{k\leq n}
Z_k(\omega)^2 \Bigr] +
3\mathbb{E}\Bigl[ \sup_{k\leq n} R_k^2 \Bigr].
\]
Since $M_n$ is a martingale, Doob's inequality and the first statement
in \eqref{MnZnL2} imply that
\[
\mathbb{E}\Bigl[ \sup_{k\leq n} M_k^2 \Bigr] \leq4 \mathbb
{E}[ M_n^2 ] =
\mathcal{O}(n).
\]
The same argument given above which showed that $\mathbb{E}R_n^2 =
o(n)$ can
be repeated to show that, for any $\beta\in(1/2,1)$, there exists a
constant $C<\infty$ such that
\[
\mathbb{E}\Bigl[ \sup_{k\leq n} R_k^2 \Bigr] \leq C n^{\beta} + C
n \mathbb{P}\Bigl( \sup
_{k\leq n} |X_k-k\mathrm{v}_P| \geq n^{\beta} \Bigr) = o(n),
\]
where in the last equality we used the averaged functional central
limit theorem.
% The averaged functional central limit theorem follows from the
%quenched functional central limit theorem, the convergence of $Z_{nt}(
%$E_\w T_1$.
% Instead of using averaged functional CLT we could also use seperately
%the quenched functional CLT and the convergence of $Z_{n\cdot}(\w)$ to
%Brownian motion to prove the above probability vanishes.
To finish the proof of \eqref{supUIclaim}, we need to show that $E_P
[ \sup_{k\leq n} Z_k(\omega)^2 ] = \mathcal{O}(n)$.
Since $Z_n(\omega) = h(\lfloor n\mathrm{v}_P \rfloor,\omega)$, this
is equivalent to showing that\break
$E_P [ \sup_{k\leq n} h(k,\omega)^2 ] = \mathcal
{O}(n)$. However, H\"
older's inequality and \eqref{hLpbound} imply that there exists an
$\eta
>0$ and $C<\infty$ such that
\[
E_P \Bigl[ \sup_{k\leq n} h(k,\omega)^2 \Bigr] \leq\Bigl( E_P
\Bigl[ \sup
_{k\leq n} |h(k,\omega)|^{2+2\eta} \Bigr] \Bigr)^{1/(1+\eta)}
\leq Cn.
\]
This completes the proof of Proposition \ref{UIprop}.
\end{appendix}

\printaddresses


\begin{thebibliography}{23}

%b1 ###
\bibitem{arratia}
\begin{barticle}[mr]
\bauthor{\bsnm{Arratia},~\bfnm{Richard}\binits{R.}}
(\byear{1983}).
\btitle{The motion of a tagged particle in the simple symmetric exclusion
  system on {$\textbf{Z}$}}.
\bjournal{Ann. Probab.}
\bvolume{11}
\bpages{362--373}.
\bid{mr={690134}}
\end{barticle}
\endbibitem

%b2 ###
\bibitem{bala-sepp-aom}
\begin{bmisc}[auto:SpringerTagBib|2009-01-14|16:51:27]
\bauthor{\bsnm{Bal{\'a}zs},~\bfnm{M.}\binits{M.}} \AND
  \bauthor{\bsnm{Sepp{\"a}l{\"a}inen},~\bfnm{T.}\binits{T.}}
  (\byear{2010}).
\bhowpublished{Order of current variance and diffusivity in the
  asymmetric simple exclusion process. \textit{Ann. of Math.} \textbf{171} 1237--1265}.
\end{bmisc}
\endbibitem

%b4 ###
\bibitem{durr-gold-lebo}
\begin{barticle}[mr]
\bauthor{\bsnm{D{\"u}rr},~\bfnm{Detlef}\binits{D.}},
  \bauthor{\bsnm{Goldstein},~\bfnm{Sheldon}\binits{S.}} \AND
  \bauthor{\bsnm{Lebowitz},~\bfnm{Joel~L.}\binits{J.~L.}}
(\byear{1985}).
\btitle{Asymptotics of particle trajectories in infinite one-dimensional
  systems with collisions}.
\bjournal{Comm. Pure Appl. Math.}
\bvolume{38}
\bpages{573--597}.
\bid{doi={10.1002/cpa.3160380508}, mr={803248}}
\end{barticle}
\endbibitem

%b3 ###
\bibitem{durr}
\begin{bbook}[mr]
\bauthor{\bsnm{Durrett},~\bfnm{Richard}\binits{R.}}
(\byear{1996}).
\btitle{Probability: Theory and Examples}, \bedition{2nd} ed.
\bpublisher{Duxbury Press}, \baddress{Belmont, CA}.
\bid{mr={1609153}}
\end{bbook}
\endbibitem



%b5 ###
\bibitem{ferr-spoh-06}
\begin{barticle}[mr]
\bauthor{\bsnm{Ferrari},~\bfnm{Patrik~L.}\binits{P.~L.}} \AND
  \bauthor{\bsnm{Spohn},~\bfnm{Herbert}\binits{H.}}
(\byear{2006}).
\btitle{Scaling limit for the space--time covariance of the stationary totally
  asymmetric simple exclusion process}.
\bjournal{Comm. Math. Phys.}
\bvolume{265}
\bpages{1--44}.
\bid{doi={10.1007/s00220-006-1549-0}, mr={2217295}}
\end{barticle}
\endbibitem

%b6 ###
\bibitem{gQCLT}
\begin{barticle}[mr]
\bauthor{\bsnm{Goldsheid},~\bfnm{Ilya~Ya.}\binits{I.~Y.}}
(\byear{2007}).
\btitle{Simple transient random walks in one-dimensional random environment:
  The central limit theorem}.
\bjournal{Probab. Theory Related Fields}
\bvolume{139}
\bpages{41--64}.
\bid{doi={10.1007/s00440-006-0038-x}, mr={2322691}}%
\end{barticle}%
\endbibitem%

%b7 ###
\bibitem{jara-09}
\begin{bmisc}[auto:SpringerTagBib|2009-01-14|16:51:27]
\bauthor{\bsnm{Jara},~\bfnm{M.}\binits{M.}}
(\byear{2009}).
\bhowpublished{Current and density fluctuations for interacting particle
systems with anomalous diffusive behavior. Available at}
\href{http://arxiv.org/abs/arXiv:0901.0229}{arXiv:0901.0229}.
\end{bmisc}
\endbibitem

%b8 ###
\bibitem{jara-landim-06}
\begin{barticle}[mr]
\bauthor{\bsnm{Jara},~\bfnm{M.~D.}\binits{M.~D.}} \AND
  \bauthor{\bsnm{Landim},~\bfnm{C.}\binits{C.}}
(\byear{2006}).
\btitle{Nonequilibrium central limit theorem for a tagged particle in symmetric
  simple exclusion}.
\bjournal{Ann. Inst. H. Poincar\'e Probab. Statist.}
\bvolume{42}
\bpages{567--577}.
\bid{doi={10.1016/j.anihpb.2005.04.007}, mr={2259975}}
\end{barticle}
\endbibitem

%b9 ###
\bibitem{jara-landim-08}
\begin{barticle}[mr]
\bauthor{\bsnm{Jara},~\bfnm{M.~D.}\binits{M.~D.}} \AND
  \bauthor{\bsnm{Landim},~\bfnm{C.}\binits{C.}}
(\byear{2008}).
\btitle{Quenched non-equilibrium central limit theorem for a tagged particle in
  the exclusion process with bond disorder}.
\bjournal{Ann. Inst. H. Poincar\'e Probab. Statist.}
\bvolume{44}
\bpages{341--361}.
\bid{doi={10.1214/07-AIHP112}, mr={2446327}}
\end{barticle}
\endbibitem

%b10 ###
\bibitem{joha}
\begin{barticle}[mr]
\bauthor{\bsnm{Johansson},~\bfnm{Kurt}\binits{K.}}
(\byear{2000}).
\btitle{Shape fluctuations and random matrices}.
\bjournal{Comm. Math. Phys.}
\bvolume{209}
\bpages{437--476}.
\bid{doi={10.1007/s002200050027}, mr={1737991}}
\end{barticle}
\endbibitem

%b11 ###
\bibitem{kksStable}
\begin{barticle}[mr]
\bauthor{\bsnm{Kesten},~\bfnm{H.}\binits{H.}},
  \bauthor{\bsnm{Kozlov},~\bfnm{M.~V.}\binits{M.~V.}} \AND
  \bauthor{\bsnm{Spitzer},~\bfnm{F.}\binits{F.}}
(\byear{1975}).
\btitle{A limit law for random walk in a random environment}.
\bjournal{Compositio Math.}
\bvolume{30}
\bpages{145--168}.
\bid{mr={0380998}}
\end{barticle}
\endbibitem

%b12 ###
\bibitem{kSTCP}
\begin{barticle}[mr]
\bauthor{\bsnm{Kumar},~\bfnm{Rohini}\binits{R.}}
(\byear{2008}).
\btitle{Space--time current process for independent random walks in one
  dimension}.
\bjournal{ALEA Lat. Am. J. Probab. Math. Stat.}
\bvolume{4}
\bpages{307--336}.
\bid{mr={2456971}}
\end{barticle}
\endbibitem

%b13 ###
\bibitem{mrzStable}
\begin{barticle}[mr]
\bauthor{\bsnm{Mayer-Wolf},~\bfnm{Eddy}\binits{E.}},
  \bauthor{\bsnm{Roitershtein},~\bfnm{Alexander}\binits{A.}} \AND
  \bauthor{\bsnm{Zeitouni},~\bfnm{Ofer}\binits{O.}}
(\byear{2004}).
\btitle{Limit theorems for one-dimensional transient random walks in {M}arkov
  environments}.
\bjournal{Ann. Inst. H. Poincar\'e Probab. Statist.}
\bvolume{40}
\bpages{635--659}.
\bid{doi={10.1016/j.anihpb.2004.01.003}, mr={2086017}}
\end{barticle}
\endbibitem

%b14 ###
\bibitem{peli-seth}
\begin{barticle}[mr]
\bauthor{\bsnm{Peligrad},~\bfnm{Magda}\binits{M.}} \AND
  \bauthor{\bsnm{Sethuraman},~\bfnm{Sunder}\binits{S.}}
(\byear{2008}).
\btitle{On fractional {B}rownian motion limits in one dimensional
  nearest-neighbor symmetric simple exclusion}.
\bjournal{ALEA Lat. Am. J. Probab. Math. Stat.}
\bvolume{4}
\bpages{245--255}.
\bid{mr={2448774}}
\end{barticle}
\endbibitem

%b15 ###
\bibitem{pThesis}
\begin{bmisc}[mr]
\bauthor{\bsnm{Peterson},~\bfnm{Jonathon}\binits{J.}}
(\byear{2008}).
\btitle{Limiting distributions and large deviations for random walks in random
  environments}.
\bhowpublished{Ph.D. thesis, Univ. Minnesota. Available at}
\href{http://arxiv.org/abs/arXiv:0810.0257v1}{arXiv:0810.0257v1}.
\end{bmisc}
\endbibitem

%b16 ###
\bibitem{p1LSL2}
\begin{barticle}[mr]
\bauthor{\bsnm{Peterson},~\bfnm{Jonathon}\binits{J.}}
(\byear{2009}).
\btitle{Quenched limits for transient, ballistic, sub-{G}aussian
  one-dimensional random walk in random environment}.
\bjournal{Ann. Inst. H. Poincar\'e Probab. Statist.}
\bvolume{45}
\bpages{685--709}.
\bid{doi={10.1214/08-AIHP149}, mr={2548499}}
\end{barticle}
\endbibitem

%b17 ###
\bibitem{psHydro}
\begin{bmisc}[mr]
\bauthor{\bsnm{Peterson},~\bfnm{Jonathon}\binits{J.}}
(\byear{2009}).
\btitle{Systems of one-dimensional random walks in a common random environment}.
\bhowpublished{Preprint. Available at}
\href{http://arxiv.org/abs/arXiv:0907.3680v1}{arXiv:0907.3680v1}.
\end{bmisc}
\endbibitem

%b18 ###
\bibitem{pzSL1}
\begin{barticle}[mr]
\bauthor{\bsnm{Peterson},~\bfnm{Jonathon}\binits{J.}} \AND
  \bauthor{\bsnm{Zeitouni},~\bfnm{Ofer}\binits{O.}}
(\byear{2009}).
\btitle{Quenched limits for transient, zero speed one-dimensional random walk
  in random environment}.
\bjournal{Ann. Probab.}
\bvolume{37}
\bpages{143--188}.
\bid{doi={10.1214/08-AOP399}, mr={2489162}}
\end{barticle}
\endbibitem

%b19 ###
\bibitem{quas-valko}
\begin{barticle}[mr]
\bauthor{\bsnm{Quastel},~\bfnm{Jeremy}\binits{J.}} \AND
  \bauthor{\bsnm{Valko},~\bfnm{Benedek}\binits{B.}}
(\byear{2007}).
\btitle{{$t\sp {1/3}$} {s}uperdiffusivity of finite-range asymmetric exclusion
  processes on {$\Bbb Z$}}.
\bjournal{Comm. Math. Phys.}
\bvolume{273}
\bpages{379--394}.
\bid{doi={10.1007/s00220-007-0242-2}, mr={2318311}}
\end{barticle}
\endbibitem

%b20 ###
\bibitem{sepp-rw}
\begin{barticle}[mr]
\bauthor{\bsnm{Sepp{\"a}l{\"a}inen},~\bfnm{Timo}\binits{T.}}
(\byear{2005}).
\btitle{Second-order fluctuations and current across characteristic for a
  one-dimensional growth model of independent random walks}.
\bjournal{Ann. Probab.}
\bvolume{33}
\bpages{759--797}.
\bid{doi={10.1214/009117904000000946}, mr={2123209}}
\end{barticle}
\endbibitem

%b21 ###
\bibitem{sRWRE}
\begin{barticle}[mr]
\bauthor{\bsnm{Solomon},~\bfnm{Fred}\binits{F.}}
(\byear{1975}).
\btitle{Random walks in a random environment}.
\bjournal{Ann. Probab.}
\bvolume{3}
\bpages{1--31}.
\bid{mr={0362503}}
\end{barticle}
\endbibitem

%b22 ###
\bibitem{walsh}
\begin{bincollection}[mr]
\bauthor{\bsnm{Walsh},~\bfnm{John~B.}\binits{J.~B.}}
(\byear{1986}).
\btitle{An introduction to stochastic partial differential equations}.
In \bbooktitle{\'{E}cole D'\'et\'e de Probabilit\'es de {S}aint-{F}lour,
  {XIV}---\textit{1984}}.
\bseries{Lecture Notes in Math.}
\bvolume{1180}
\bpages{265--439}.
\bpublisher{Springer}, \baddress{Berlin}.
\bid{mr={876085}}
\end{bincollection}
\endbibitem

%b23 ###
\bibitem{zRWRE}
\begin{bincollection}[mr]
\bauthor{\bsnm{Zeitouni},~\bfnm{Ofer}\binits{O.}}
(\byear{2004}).
\btitle{Random walks in random environment}.
In \bbooktitle{Lectures on Probability Theory and Statistics}.
\bseries{Lecture Notes in Math.}
\bvolume{1837}
\bpages{189--312}.
\bpublisher{Springer}, \baddress{Berlin}.
\bid{mr={2071631}}%
\end{bincollection}%
\endbibitem%

\end{thebibliography}
\end{document}